\documentclass[10pt,reqno,twoside]{amsart}
 \renewcommand{\i}{{\mathrm{i}}} 

\newcommand{\Q}{{\mathcal{Q}}}

\newcommand{\N}{{\mathbb{N}}}
\newcommand{\rk}{{\mathcal{R}}}

\newtheorem{theorem}{Theorem}
 
 \newtheorem{lemma}[theorem]{Lemma}
 
 \theoremstyle{remark}

\reversemarginpar 

\begin{document}

\title{Large gaps between the zeros of the Riemann zeta function}

\author{NATHAN NG}
 \address{Department of Mathematics and Statistics, 
  University of Ottawa, 
 585 King Edward
 Ave.,
 Ottawa,
 ON K1N 6N5}
 \email{nng362@science.uottawa.ca}
 \thanks{This research was funded in part by NSERC and 
 NSF FRG grant DMS 0244660}

\date{\today}

\begin{abstract}
\noindent  We show that the generalized Riemann hypothesis implies that there are infinitely many consecutive zeros of the
zeta function whose spacing is 2.9125 times larger than the average spacing.  This is deduced from the
calculation of the second moment of the Riemann zeta function multiplied by a Dirichlet polynomial averaged over
the zeros of the zeta function.
\end{abstract}

\maketitle

\section{Introduction}

If the Riemann hypothesis (RH) is true then the non-trivial zeros of the Riemann zeta function, $\zeta(s)$,
satisfy $1/2+i \gamma_{n}$ with $\gamma_{n} \in \Bbb R$.  Riemann noted that the argument principle implies that
number of zeros of $\zeta(s)$ in the box with vertices $0,1,1+iT,$ and $iT$ is
$N(T) \sim (T/2 \pi) \log \left(T/2 \pi e \right)$.
This implies that on average
$
   (\gamma_{n+1}-\gamma_{n})
   \approx 2  \pi / \log \gamma_{n}
$
and hence the average spacing of the sequence
$
   \hat{\gamma}_{n} = \gamma_{n} \log \gamma_{n}/2 \pi
$
is one.  Montgomery \cite{Mo} 
investigated the pair correlation of these numbers and he proposed the
fundamental conjecture
\begin{equation}
  \frac{1}{N}  \# \{ 1 \le j \ne k \le N \ | \
                  a \le \hat{\gamma}_{j} - \hat{\gamma}_{k} \le b \ \} \sim
   \int_{a}^{b} \left( 1 - \left( \frac{\sin \pi x}{\pi x} \right)^{2} \right) \, dx
   \label{eq:pairc}
\end{equation}
for $0 < a < b$ as $N \to \infty$. Moreover,
it is expected that the consecutive spacings,
$\hat{\gamma}_{n+1} - \hat{\gamma}_{n}$, have a limiting
distribution 
function which agrees with the Gaussian Unitary Ensemble from random
matrix theory.  See Odlyzko \cite{Od}
for extensive numerical evidence in favour of this conjecture 
and also see Rudnick-Sarnak \cite{RS} for a study
of the $n$-level correlations of $\hat{\gamma}_n$. 
In light of the expected distribution of the consecutive spacings of zeta
Montgomery suggested in \cite{Mo} 
that there exist arbitrarily large and small gaps
between the zeros of the zeta function.  That is to say
\[
   \lambda = \limsup_{n \to \infty} (\hat{\gamma}_{n+1} - \hat{\gamma}_{n})
   = \infty \ \mathrm{and} \
   \mu = \liminf_{n \to \infty} (\hat{\gamma}_{n+1} - \hat{\gamma}_{n})
   = 0 \ .
\]

In this article, we focus on the large gaps and we assume the generalized Riemann hypothesis (GRH) is true.  This
conjecture states that the non-trivial zeros of the Dirichlet $L$-functions are on the 
$\mathrm{Re}(s)=1/2$ line.  We
establish
\newtheorem{gaps}{{\bf Theorem}}
\begin{gaps}
The generalized Riemann hypothesis implies $\lambda
> 2.9125$.
\end{gaps}
Selberg was the first to establish that $\lambda > 1$ based on his work
concerning moments of $S(t)=(1/\pi ) \arg \zeta(1/2+it)$ in short 
intervals.  Montgomery and Odlyzko \cite{MO} obtained 
$\lambda > 1.9799$ assuming the Riemann hypothesis.  
The current record due to Hall is $\lambda > 2.34$.  
Hall's work makes use of Wirtinger's inequality in conjunction with
asymptotic formulae for continuous mixed
moments of the zeta function and its derivatives.  Moreover, Hall is currently attempting to show that the
asymptotic evaluation of all mixed moments of zeta and its derivatives yields $\lambda=\infty$.  It should be
noted that the best published result \cite{CGG0} assuming the Riemann hypothesis is worse than Hall's unconditional work.  Theorem 1 extends earlier results of Conrey, Ghosh, and Gonek where they assume GRH to obtain $\lambda > 2.68$.  In fact, their work is based on the following idea of J. Mueller\cite{Mu}.  
Let $H: \mathbb{C} \to \mathbb{R}$ and consider the associated functions
\begin{align}
  &  \mathcal{M}_{1}(H,T) = \int_{1}^{T} H(1/2+\i t) \, dt \ ,  \\
  &  m(H,T;\alpha) = \sum_{T < \gamma < 2T}
     H(1/2+ \i (\gamma + \alpha))  \ , \\
  & \mathcal{M}_{2}(H,T;c) = \int_{-c/L}^{c/L} m(H,T;\alpha) \, d
  \alpha 
\end{align}
where we put $L = \log (T/2 \pi)$.
This notation shall be used throughout the article. 
However, one notes that
\begin{equation}
  \frac{\mathcal{M}_{2}(H,2T;c) - \mathcal{M}_{2}(H,T;c)}{\mathcal{M}_{1}(H,2T) -
  \mathcal{M}_{1}(H,T)} < 1
  \label{eq:mc2}
\end{equation}
implies $\lambda > \frac{c}{\pi}$.  Mueller applied this idea with $H(s)=|\zeta(s)|^2$ and obtained $\lambda > 1.9$. 
We should note that the method of Montgomery and Odlyzko \cite{MO}
is equivalent to the method of Mueller \cite{Mu}.  This was realized
later by the authors of \cite{CGG0}.  Now consider the
Dirichlet polynomial
\begin{equation}
   A(s) = \sum_{n \le y} a(n)n^{-s}  \ .
\end{equation}
Assuming the Riemann hypothesis, Conrey, Ghosh, and Gonek in \cite{CGG0} applied~(\ref{eq:mc2}) to $H(s)=|A(s)|^2$ with $a(n)=d_{2.2}(n)$,
$y=T^{1-\epsilon}$ and obtained $\lambda > 2.337$ (and $\mu < 0.5172$).
Here $d_{r}(n)$ is the coefficient of $n^{-s}$ in the Dirichlet series
$\zeta(s)^r$.  If $r$
is a natural number then $d_{r}(n)$ equals the number of representations of $n$ as a product of $r$ positive
integers.  In recent work \cite{Ng}, we have shown that the Riemann hypothesis implies $\lambda > 2.56$ (and $\mu < 0.5162$). 
In \cite{CGG}, Conrey, Ghosh, and Gonek applied~(\ref{eq:mc2}) to $H(s)=|\zeta(s)A(s)|^2$ with $a(n)=1$ and $y=(T/2
\pi)^{\frac{1}{2}-\epsilon}$ and obtained $\lambda > 2.68$.
However, in this situation it is necessary to assume GRH 
in order to evaluate the discrete mean value $m(H,T;\alpha)$.  
We continue this programme by considering a more general choice for the
coefficient $a(n)$. Precisely, we choose as our function 
$H_{r}(s)=|\zeta(s)A(s)|^{2}$
where $A(s)$ has coefficients
\begin{equation}
   a(n) = d_{r}(n) P \left( \frac{\log n}{\log y} \right)\
\end{equation}
for $P$ a polynomial and for $r \in \mathbb{N}$.  Furthermore, we put
$y = \left( T/2 \pi\right)^{\eta}$
where $\eta < 1/2$.  Ideally, we would like to evaluate $m(|\zeta(s)A(s)|^2,T;\alpha)$ with arbitrary
coefficients.  Our calculation follows that of \cite{CGG}.  However,
we must take into account that $d_{r}$
is not a completely multiplicative function for $r \ge 2$.  It should be noted that Chris Hughes \cite{Hu} has shown that if
$H(s)=|\zeta(s)|^{4}$ is admissible then his random matrix theory conjectures yield $\lambda > 2.7$.  In addition,
he has shown (unpublished) that if $H(s) = |\zeta(s)|^{k}$ is admissible for arbitrarily large $k$ then the random matrix theory
conjectures for $\mathcal{M}_1(H,T)$ and $m(H,T;\alpha)$ yield $\lambda \ge f(k)$
where $f(k) \nearrow \infty$ at a linear rate.  By choosing
$H_{r}(s)=|\zeta(s)A_{r}(s)|^2$ with coefficients
$a(n)=d_{r}(n)$ we are hoping that $H_{r}(s)$ will mimic the larger moment $|\zeta(s)|^{2r+2}$.  The work
of \cite{CGG} corresponds to the choice $r=1,P(x)=1$.

We now state the precise result.  We define several functions that will appear in the course of the proof. Given a
polynomial $P$ and $u \in \mathbb{Z}_{\ge 0}$ we define
\begin{equation}
   Q_{u}(x) = \int_{0}^{1} \theta^{u} P(x + \theta(1-x)) \, d \theta  \ .
\end{equation}
Given $\vec{n} = (n_{1},n_{2},n_{3},n_{4},n_{5}) \in (\mathbb{Z}_{\ge 0})^{5}$ we define
\begin{equation}
  i_{P}(\vec{n}) = \int_{0}^{1} \int_{0}^{1-x}
  x^{r^2-1}(1-x)^{n_{1}} (1-y-x)^{n_{2}} y^{n_{3}} Q_{n_{4}}(x) Q_{n_{5}}(x+y) dy dx  \ .
  \label{eq:i1234}
\end{equation}
For $\eta \in \mathbb{R}$ and $\vec{n}=(n_{1},n_{2},n_{3})\in (\mathbb{Z}_{\ge 0})^{3}$ we define
\begin{equation}
  k_{P}(\vec{n}) = k_{P}(\vec{n};\eta) = \int_{0}^{1} \int_{0}^{1-x} x^{r-1}(\eta^{-1}-x)^{n_{1}}
  y^{r^2-1} (1-y)^{n_{2}}P(x+y) Q_{n_{3}}(y) dy dx \ .
  \label{eq:k1234}
\end{equation}
Recall $\eta$ corresponds to the length of our Dirichlet polynomial.  Given $r \ge 1$ we define the constants
\begin{equation}
  a_{r} = \prod_{p} \left(
  (1-p^{-1})^{r^2} \sum_{m=0}^{\infty} \left(
  \frac{\Gamma(r+m)}{\Gamma(r) m!}
  \right)^{2} p^{-m}
  \right)   \  \mathrm{and} \
   C_{r} =  \frac{a_{r+1}}{(r^2-1)!((r-1)!)^2} \ . 
  \label{eq:ar}
\end{equation}
With all of these definitions in hand we present our result for $m(H_{r},T;\alpha)$.
\newtheorem{dismom}[gaps]{{\bf Theorem}}
\begin{dismom}
Suppose $r \in \mathbb{N}$ and $\eta < 1/2$. GRH implies
\begin{equation}
\begin{split}
  & m(H_{r},T;\alpha) \sim \frac{C_{r}TL^{(r+1)^2+1}}{\pi} 
   \mathrm{Re}
  \sum_{j=1}^{\infty}  z^{j}
  \eta^{j+(r+1)^2+1} \left(  \frac{ r \hat{i}(r,j,\eta)}{j!} +
   \hat{k}(r,j,\eta) \right)
  \label{eq:thm2}
\end{split}
\end{equation}
where $z = i \alpha L$, $|z| \ll 1$,
\begin{equation}
   \hat{i}(r,j,\eta) = -i_{P}(r,r,j,r-1,r-1)\eta^{-1}+ i_{P}(r+1,r,j,r,r-1)+i_{P}(r,r+1,j,r-1,r)
    \ ,
\end{equation}
\begin{equation}
\begin{split}
  \hat{k}(r,j,\eta) & = -(r-1)!\sum_{n=-2}^{\min (j,r-2)}
  \frac{(-1)^{n} \binom{r}{n+2}}{(j-n)! (r+n+1)!}
  k_{P}(j-n,r+n+2,r+n+1)   \ .
\end{split}
\end{equation}
This result is valid up to an error term which is 
$O_{\epsilon,r}(TL^{(r+1)^2}+T^{1/2+\eta+\epsilon})$.
\end{dismom}
We note that it is probable that Theorem 2 can be proven only assuming
the Generalized Lindel\"{o}f Hypothesis by following the work of Conrey,
Ghosh, and Gonek \cite{CGG2} on simple zeros of $\zeta(s)$.  
Even this assumption may possibly be weakened further since the
main theorem in \cite{CGG2}
actually assumes an upper bound for the sixth integral moment of $L(s,\chi)$ 
on average. 
Also we remark that the case $r=1,P(x)=1$ reduces, after some calculation, to
\begin{equation}
\begin{split}
  m(H_{1},T;\alpha) & \sim 
  \frac{6}{\pi^{2}} \frac{TL^{5}}{2 \pi} 
  \sum_{j=0}^{\infty} \frac{(-1)^{j+1}(\alpha L)^{2j+2}}{(2j+5)!} \cdot  \\
  & \cdot 
  \left(
  \frac{-3 \eta^{2} + (2j+5) \eta^{3}}{3} - \frac{2j+5}{j+3} \eta^{2j+6}
  + \eta^{2j+7} + \eta^{2}(1-\eta)^{2j+5} 
  \right)   \\
\end{split}
\end{equation}
which corresponds to Theorem 1 of \cite{CGG}.   \\
\\
\noindent {\it Acknowledgements}  The research for this article commenced
at the Universit\'e  de Montr\'eal during the 2003-2004
academic year and was completed the following year at the University of Michigan.  The author thanks Chris Hughes for 
correspondence concerning his unpublished work.

\section{Theorem 2 implies Theorem 1}

In this section, we deduce Theorem 1 from Theorem 2.  The rest of
the article will be devoted to establishing the discrete moment result
of Theorem 1.  Put $\eta=1/2-\epsilon$ with $\epsilon$ arbitrarily small. Since $\mathrm{Re}(z^{j}) = (-1)^{k}
(\alpha L)^{2k}$ if $j=2k$ and zero otherwise, it follows from~(\ref{eq:thm2}) that
\begin{equation}
  m(H_{r},2T;\alpha) -m(H_{r},T;\alpha) = \phi(r,\eta, \alpha) \frac{C_{r}TL^{(r+1)^2+1}}{\pi}
  (1+O(L^{-1}))
  \label{eq:m12}
\end{equation}
where
\[
  \phi(r,\eta,\alpha) = \eta^{(r+1)^2+1}
  \sum_{j=1}^{\infty} (-1)^{j}  (\alpha L \eta)^{2j}
  \left( \frac{r\hat{i}(r,2j,\eta)}{(2j+1)!}
  + \frac{\hat{k}(r,2j,\eta)}{2j+1} \right)
 \ .
\]
Integrating~(\ref{eq:m12}) with respect to $\alpha$ over the interval $[-c/L,c/L]$ we have
$\mathcal{M}_{2}(H_r,2T;c) - \mathcal{M}_{2}(H_r,T;c)$ equals
\[
  \frac{2C_{r}TL^{(r+1)^2}\eta^{(r+1)^2+1}}{\pi}
  \sum_{j=1}^{\infty} (-1)^{j}c^{2j+1} \eta^{2j} \left( \frac{r\hat{i}(r,2j,\eta)}{(2j+1)!}
  + \frac{\hat{k}(r,2j,\eta)}{2j+1} \right)
\]
plus an error $O(TL^{(r+1)^2})$.  In the above expression, we may replace $\eta=1/2-\epsilon$
by $1/2$ yielding
\begin{equation}
\begin{split}
 & \mathcal{M}_{2}(H_{r},2T;c) - \mathcal{M}_{2}(H_{r},T;c) = \frac{2C_{r}TL^{(r+1)^2} \eta^{(r+1)^2+1}}{\pi} \cdot \\
 & \sum_{j=1}^{\infty} \frac{(-1)^{j}c^{2j+1}}{2^{2j}}
 \left( \frac{r\hat{i}(r,2j,\frac{1}{2})}{(2j+1)!}
  + \frac{\hat{k}(r,2j,\frac{1}{2})}{(2j+1)}
  \right)
  + O(\epsilon TL^{(r+1)^2}) \ .
  \nonumber
\end{split}
\end{equation}
We now recall the following result of Conrey and Ghosh \cite{CG}.
\begin{lemma} If $y=T^{\eta}$ with $0 < \eta < 1/2$ then
\begin{equation}
\begin{split}
    & \mathcal{M}_{1}(H_r,T)   \sim 
    \frac{a_{r+1}}{((r-1)!)^{2}(r^{2}-1)!} T(\log y)^{(r+1)^{2}}
     \cdot  \\
     &  \int_{0}^{1}\alpha^{r^2-1}
     (\eta^{-1}  (1-\alpha)^{2r} Q_{r-1}(\alpha)^2
     -2  (1-\alpha)^{2r+1} Q_{r}(\alpha)Q_{r-1}(\alpha))
     \, d \alpha  
     \label{eq:M1Hr}
\end{split}
\end{equation}
as $T \to \infty$.  This is valid up to an error term which is $O(L^{-1})$
smaller than the main term. 
\end{lemma}
\noindent Hence, we have
\begin{equation}
\begin{split}
  & \mathcal{M}_{1}(H_{r},2T) - \mathcal{M}_{1}(H_{r},T) = 
  C_{r} T(L\eta)^{(r+1)^2} \\
  & \int_{0}^{1} \alpha^{r^2-1} ( \eta^{-1}(1-\alpha)^{2r}Q_{r-1}(\alpha)^{2}
   -2(1-\alpha)^{2r+1} Q_{r-1}(\alpha)Q_{r}(\alpha)) \, d \alpha
  + O(\epsilon TL^{(r+1)^2} ) \ .
  \nonumber
\end{split}
\end{equation}
We deduce that
\[
 \frac{\mathcal{M}_{2}(H_{r},2T;c)-\mathcal{M}_{2}(H_{r},T;c)}{\mathcal{M}_{1}(H_{r},2T)
  -\mathcal{M}_{1}(H_{r},T)} = f_{r}(c) + O(\epsilon)
\]
where
\begin{align}
  f_{r}(c) & = \frac{1}{D}     
  \sum_{j=1}^{\infty} \frac{(-1)^{j} c^{2j+1}}{2^{2j}}
  \left( \frac{r\hat{i}(r,2j,\frac{1}{2})}{(2j+1)!}
  + \frac{\hat{k}(r,2j,\frac{1}{2})}{2j+1} \right) 
  \label{eq:frc}
\end{align}
and 
\[
    D :=  \pi \int_{0}^{1} \alpha^{r^2-1} ( \eta^{-1}(1-\alpha)^{2r}Q_{r-1}(\alpha)^{2}
   -2(1-\alpha)^{2r+1} Q_{r-1}(\alpha)Q_{r}(\alpha)) \, d \alpha
\]
We define $\lambda_{r} := \sup_{f_{r}(c) < 1} ( c )$
and thus $\lambda \ge \frac{\lambda_{r}}{\pi}$.
We may now compute~(\ref{eq:frc}) for various choices of $r$ and $P(x)$.
For example,  we shall choose $c=2.9125 \pi$, $r=2$ and $P(x)=1-0.1x+100x^2
-0.2x^3$.
We compute the sum as follows: by a Maple calculation we have 
\[
    D^{-1} \sum_{j=0}^{J}\frac{(-1)^{j} c^{2j+1}}{2^{2j}}
  \left( \frac{2\hat{i}(2,2j,\frac{1}{2})}{(2j+1)!}
  + \frac{\hat{k}(2,2j,\frac{1}{2})}{2j+1} \right)  =
  0.9999845837
\]
for $J=80$.  On the other hand, we may bound the terms $j>J$.  
Since $|Q_{u}(x)| \le ||P||_{1}$ we may establish the crude bound
\[
   |i_{P}(\vec{n})| \le 
   \frac{||P||_{1}^{2} (r^2-1)! (n_1+n_3+1)!}{(n_1+n_3+r^2+1)!(n_3+1)} 
\]    
for $\vec{n} \in ( \mathbb{Z}_{\ge 0})^{5}$.
It thus follows that 
\[
  |\hat{i}(r,2j,1/2)|
  \le \frac{||P||_1^2 (r^2-1)!}{2j+1}
  \left(
  \frac{4 (r+2j+1)!}{(r^2+r+2j+1)!}
  \right)
\]
and hence
\begin{equation}
\begin{split}
  & \left| \frac{1}{D}     
  \sum_{j>J}^{\infty} \frac{(-1)^{j} c^{2j+1}}{2^{2j}}
  \frac{2\hat{i}(2,2j,\frac{1}{2})}{(2j+1)!}    \right| \le  
  \frac{48c||P||_1^2}{D(2J)} \sum_{j > J} 
   \frac{(c/2)^{2j}(2j+3)!}{(2j+1)!(2j+7)!}  \\
   & \le \frac{48c ||P||_{1}^2}{\sqrt{2 \pi} D (2J)^5}
   \sum_{j > J} e^{- 2j ( \log(2j)-(\log(c/2)+1))} \\
   & <   \frac{48c ||P||_{1}^2}{\sqrt{2 \pi} D (2J)^5} 
   \frac{e^{-2J(\log(2J)-\log(c/2)-1)}}{2(\log(2J)-\log(c/2)-1)}  < 10^{-184}
   \nonumber
\end{split}
\end{equation}
where we have applied $n! > (n/e)^{n}$.  
A similar calculation establishes that
\[
    \left| \frac{1}{D}     
  \sum_{j>J}^{\infty} \frac{(-1)^{j} c^{2j+1}}{2^{2j}}
  \frac{\hat{k}(2,2j,\frac{1}{2})}{(2j+1)}    \right| <  10^{-130} \ . 
\]
We conclude that $f_{2}(2.1925 \pi) < 1$ and hence establish Theorem 1.
We made our choice of $r$ and $P(x)$ by a computer search. 
We note that there are many choices of $r$ and $P(x)$ that improve
the work of \cite{CGG}.  For example, $r=3$, $P(x)=1$ yields
$\lambda > 2.78$ whereas $r=2$, $P(x)=1$ yields $\lambda > 2.86$.

\section{Some notation and definitions}

Throughout this article we shall employ that notation 
\begin{equation}
    [t]_{y} := \frac{\log t}{\log y}  
\end{equation}
for $t,y>0$.  This will allow us to write several
equations more compactly.  In addition, we shall encounter a
variety of arithmetic functions.  We define $j_{\tau}(n)$,
$\Lambda(n)$, and $d_{r}(n)$ as follows:   
\begin{equation}
   j_{\tau}(n) = \prod_{p \mid n} (1+O(p^{-\tau})) 
   \label{eq:jtau}
\end{equation}
for $\tau >0$ and the constant in the $O$ is fixed and independent of $\tau$. 
Next $\Lambda(n)$ and $d_{r}(n)$ may be defined by their Dirichlet 
series generating functions:
\[
    - \frac{\zeta^{'}(s)}{\zeta(s)} = \sum_{n=1}^{\infty}\frac{\Lambda(n)}{n^s}
    \ \mathrm{and} \
    \zeta(s)^{r} = \sum_{n=1}^{\infty} \frac{d_{r}(n)}{n^s} \ . 
\]
Since this article concerns the calculation of 
discrete mean values of $m(H_r,T,\alpha)$ we need to invoke several
properties of $d_r$.  Throughout this article we apply 
repeatedly the following facts concerning $d_r$:
\begin{equation}
\begin{split}
    & \sum_{a =0}^{\infty} d_{r}(p^{a})p^{-as} =
   \left(1- p^{-s} \right)^{-r}  \ ,  \\
   &  \sum_{m \le x} d_{r}(m) m^{-1} \ll \log^{r} x \ ,  \ \mathrm{and} \\
    & \sum_{m \le x}d_{r}(m)^{2} m^{-1} \ll \log^{r^2} x \ . 
    \label{eq:ram}
\end{split}
\end{equation}
In hindsight, we realize that there is nothing really special about the 
multiplicative function $d_r$ and that the calculation of this article
can be a done for more general multiplicative functions $f$ subject to
 certain simple assumptions.

\section{Initial manipulations}
In this section we set up the plan of attack for our evaluation of $m(H_{r},T;\alpha)$.  Recall that $T$ is large,
$L = \log(T/2 \pi)$, and $\epsilon$ can be made arbitrarily small.
Let $R$ denote the positively oriented contour with vertices $a+i,a+i(T+\alpha)$, $1-a+i(T+\alpha)$, $1-a+i$, the top edge
of which has a small semicircular indentation centred at $1/2+i(T+\alpha)$ opening downward and $a=1+O(L^{-1})$.
By an application of Cauchy's residue theorem, the reflection principle, and RH we have
\[
  m(H_{r},T;\alpha) = \frac{1}{2 \pi i} \int_{R}
  \frac{\zeta^{'}}{\zeta}(s-i\alpha) \zeta(s) \zeta(1-s)A(s)A(1-s)
  \, ds \ .
\]
For $s$ in the interior or boundary of $R$ we have 
$A(s) \ll_{\epsilon} y^{1-\sigma+\epsilon}$ and
 $\zeta(s)  
    \ll_{\epsilon} T^{1/2(1-\sigma) + \epsilon}$. 
The first bound is elementary and the second is the convexity bound. These combine
to give 
$\zeta(s) \zeta(1-s) A(s) A(1-s) \ll y T^{1/2+\epsilon}$.
Now choose $T^{'}$ such that $T-2 < T^{'} < T-1$ such that $T^{'}+\alpha$ is not
the ordinate of a zero of $\zeta(s)$ and 
$ (\zeta^{'}/\zeta)(\sigma+iT^{'}) \ll L^{2}$,
uniformly for $-1 \le \sigma \le 2$.  A simple argument using Cauchy's residue theorem establishes that the top edge of the contour is 
$yT^{1/2+\epsilon}$.  Similarly,  the bottom edge of the contour is $\ll_{\epsilon} yT^{\epsilon}$ since $|\zeta(s)| \ll 1$  
for $|s| \ll 1$ and $|s-1| \gg 1$.  Differentiating the functional equation,
$\zeta(1-s) = \chi(1-s)\zeta(s)$, we have
\begin{equation}
   \frac{\zeta^{'}}{\zeta}(1-s-i \alpha) = \frac{\chi^{'}}{\chi}(1-s-i \alpha) - \frac{\zeta^{'}}{\zeta}(s+i \alpha) \ .
   \label{eq:dfe} 
\end{equation}
where $\chi(s) = 2^{s} \pi^{s-1} \sin (\pi s/2) \Gamma(1-s)$.
Now the right edge is 
\begin{equation}
   I = \frac{1}{2 \pi i} \int_{a+i}^{a+i(T+\alpha)}
    \frac{\zeta^{'}}{\zeta}(s-i\alpha) \zeta(s) \zeta(1-s)A(s)A(1-s)
  \, ds
  \label{eq:I}
\end{equation}
and the left edge is by~(\ref{eq:dfe})  
\begin{equation}
\begin{split}
  & \frac{1}{2 \pi i} \int_{1-a+i(T+\alpha)}^{1-a+i}
    \frac{\zeta^{'}}{\zeta}(s-i\alpha) \zeta(s) \zeta(1-s)A(s)A(1-s)
  \, ds
 \\
  & = \frac{1}{2 \pi i}
  \int_{a-i(T+\alpha)}^{a-i}
   \left( \frac{\zeta^{'}}{\zeta}(s+i\alpha) - \frac{\chi^{'}}{\chi}(1-s-i\alpha)   
  \right) \zeta(s) \zeta(1-s)A(s)A(1-s)
  \, ds
 \\
  & =  \overline{I}-\overline{J} 
  \nonumber
\end{split}
\end{equation}
\begin{equation}
    \mathrm{where} \ 
   J = \frac{1}{2 \pi i} \int_{a+i}^{a+i(T+\alpha)}
    \frac{\chi^{'}}{\chi}(1-s+i\alpha) \zeta(s) \zeta(1-s)A(s)A(1-s)
  \, ds \ .
  \label{eq:J}
\end{equation}
Combining results we obtain
\begin{equation}
   m(H_{r},T;\alpha)=2 \mathrm{Re} I - \overline{J} +  O_{\epsilon}(yT^{\frac{1}{2}+\epsilon}) \ .
   \label{eq:mHr}
\end{equation}
We begin with the evaluation of $J$ since it is rather simple.  By Stirling's formula one has $(\chi^{'}/\chi)(1-s+i\alpha)= -\log(t/2 \pi) + O(t^{-1})$
for $t \ge 1$, $1/2 \le \sigma \le 2$, and $|\alpha| \le cL^{-1}$. 
By moving the contour to the $1/2$ line in~(\ref{eq:J}) and 
then substituting the previous estimate we obtain
\[
   J = - \frac{1}{2 \pi}
     \int_{1}^{T} ( \log t/(2 \pi)) |\zeta A(1/2+it)|^{2} \, dt
     + O \left(
      \int_{1}^{T} |\zeta A(1/2+it)|^{2} \, \frac{dt}{t}
     +yT^{\frac{1}{2}+\epsilon} \right) \ .
\]
The last term comes from the horizontal integral. An integration by parts shows that the second integral is
$L^{(r+1)^2+1}$ and therefore
$$
   J=-\frac{L}{2 \pi} \mathcal{M}_{1}(H_r,T)
   + \int_{1}^{T} \mathcal{M}_{1}(H_r,t) \frac{dt}{t}
   +O(L^{(r+1)^2+1}+yT^{\frac{1}{2}+\epsilon}))
$$
where 
\[
    \mathcal{M}_{1}(H_r,T) = 
    \int_{1}^{T} H_{r}(1/2+it) \, dt = 
    \int_{1}^{T} |\zeta(1/2+it)|^{2}
               |A(1/2+it)|^{2} \, dt
\]
By Lemma 1 above, we thus deduce 
\begin{equation}
\begin{split}
   J  \sim - \frac{C_{r}TL^{(r+1)^2+1}}{2 \pi} &
   \left(  \eta^{(r+1)^2-1}
   \int_{0}^{1} \alpha^{r^2-1}(1-\alpha)^{2r} Q_{r-1}(\alpha)^2 \, d \alpha
   \right. \\
   & \left. - 2 \eta^{(r+1)^2}
     \int_{0}^{1} \alpha^{r^2-1}(1-\alpha)^{2r+1}
      Q_{r-1}(\alpha) Q_{r}(\alpha) \, d \alpha \right)
      \label{eq:J2}
\end{split}
\end{equation}
which is valid up to an error term $O(L^{-1})$ smaller.
We have now reduced the evaluation of $m(H_{r},T;\alpha)$ to that of $I$.
We begin our evaluation of $I$ with some intial simplifications.  By the functional equation~(\ref{eq:I})
becomes
\[
   I = \frac{1}{2 \pi  i}  \int_{a+i}^{a+i(T+\alpha)}
         \chi(1-s)B(s)A(1-s) \, ds
\]
where
$B(s) =  \frac{\zeta^{'}}{\zeta} (s-i\alpha) \zeta^{2}(s) A(s) =
   \sum_{j=1}^{\infty} b(j) j^{-s}$
and
\begin{equation}
   b(j) = -   \sum_{{\begin{substack}{hmn=j
         \\ h \le y}\end{substack}}} d_{r}(h) P \left( [h]_{y} \right)
    d(m) \Lambda(n)n^{i\alpha}  \ .
   \label{eq:bj}
\end{equation}
However, Lemma 2 pp. 504-506 of \cite{CGG} deals with such integrals. 
\begin{lemma} Suppose $B(s)=\sum_{j \ge 1}b(j)j^{-s}$ and
$A(s)=\sum_{k \le y} a(k)k^{-s}$ where $a(j) \ll d_{r_1}(j)(\log j)^{l_1}$
and $b(j) \ll d_{r_2}(j) (\log j)^{l_2}$ for some non-negative integers $r_1,
r_2,l_1,l_2$
and $T^{\epsilon} \ll y \ll T$ for some $\epsilon > 0$. If
$$
   I = \int_{c+i}^{c+iT} \chi(1-s)B(s)A(1-s) \, ds
$$
then
$$
   I = \sum_{k \le y} \frac{a(k)}{k}
       \sum_{j \le \frac{nT}{2 \pi}} b(j) e \left(- j/k \right)
       + O (y T^{\frac{1}{2}} (\log T)^{r_1+r_2+l_1+l_2})  \ . 
$$
\end{lemma}
\noindent We deduce that
\begin{equation}
  I = \sum_{k \le y} \frac{d_{r}(k)P([k]_{y})}{k}
  \sum_{j \le \frac{kT}{2 \pi}} b(j) e \left(- j/k \right)
       + O (y T^{\frac{1}{2}+\epsilon})  \ . 
  \label{eq:Iexp}
\end{equation}

The goal of the rest of this paper is to evaluate the sum in~(\ref{eq:Iexp}).  
We now give a brief sketch how the proof shall proceed.  We define the Dirichlet series
\[
      Q^{*}(s,\alpha,k) = \sum_{j=1}^{\infty} b(j) e(-j/k) j^{-s} \ . 
\]  
Now the inner sum in~(\ref{eq:Iexp}) can be written by Perron's formula
as 
\begin{equation}
     \frac{1}{2 \pi i} \int_{(c)}  Q^{*}(s,\alpha,k) \left( \frac{kT}{2 \pi} \right)^{s} \frac{ds}{s} 
     = M(k) + E(k)
     \label{eq:percon}
\end{equation}
with $c > 1$.   We shall move this contour left to $\mathrm{Re}(s)=1/2+L^{-1}$
and we will have a main term, $M(k)$, arising from the residues of 
$Q^{*}(s,\alpha,k)$ at $s=1,s=1+i\alpha$.   Morover, the contribution from 
the contour on the line $\mathrm{Re}(s)=1/2+L^{-1}$ will be an error term
denoted $E(k)$.  Next $M(k)$ will be reinserted in~(\ref{eq:Iexp}) and this
will give the main term in the evaluation of $I$.  The rest of the proof 
concentrates on the calculation of 
\begin{equation}
      I =  \sum_{k \le y} \frac{d_{r}(k)P([k]_{y})}{k} M(k)
      \label{eq:Iexp2}
\end{equation}
and this part of the calcuation will be somewhat complicated. However,
it should be noted that the evaluation of~(\ref{eq:Iexp2}) will not require
GRH as it is essentially an elementary arithmetic sum.  

We now explain the connection to the Generalized Riemann Hypothesis
and how it will be invoked in the argument.   Note that the 
additive character $e(-j/k)$ may be written in terms of multiplicative characters.  In particular, if 
$(j,k)=1$ we have
the nice formula
\begin{equation}
  e(-j/k) = \frac{1}{\phi(k)} \sum_{\chi (\mathrm{mod} \ k)}
  \overline{\chi}(-j) \tau(\chi)  \ . 
  \label{eq:ejk}
\end{equation}
By this identity we shall decompose $Q^{*}(s,\alpha,k)$ into combinations of 
$L(s,\chi)$ and its logarithmic derivative where $\chi$ is a character mod
$l$ for $l \mid k$.  Now by assuming GRH we guarantee that $Q^{*}(s,\alpha,k)$ 
has only the poles at $s=1,1+i\alpha$.  If GRH were false then there would be 
extra poles occuring at those zeros that violate GRH.  This obviously would complicate the argument.  Secondly,  we require a Lindel\"{o}f type bound for $L(s,\chi)$ and 
$(L^{'}/L)(s,\chi)$ in order to ensure that the error term $E(k)$ in~(\ref{eq:percon}) is small.  Finally, we mention that many of the technicalities in 
evaluating~(\ref{eq:Iexp2}) arise from the fact that~(\ref{eq:ejk}) only
holds for $(j,k)=1$.

\section{Lemmas}

In this section we present the lemmas that will be required for 
the bounding the contribution coming from the error terms, $E(k)$,  and for evaluating the main term~(\ref{eq:Iexp2}). 
The next lemma is useful for analyzing Dirichlet series that
are products of several other Dirichlet series. 
\begin{lemma} Suppose that $
   A_{j}(s) = \sum_{n=1}^{\infty} \alpha_{j}(n) n^{-s}
$ is absolutely convergent for $\sigma  > 1$, for $1 \le j \le J$,
and that
$$
   A(s) = \sum_{n =1}^{\infty} \frac{\alpha(n)}{n^{s}}
   = \prod_{j=1}^{J} A_{j}(s) \ .
$$
Then for any positive integer $d$,
$$
   \sum_{n =1}^{\infty} \frac{\alpha(dn)}{n^{s}}
   = \sum_{d_{1} \cdots d_{J}=d}
   \prod_{j=1}^{J} \left(
    \sum_{{\begin{substack}{n=1
         \\ (n,P_{j})=1}\end{substack}}}^{\infty}
    \frac{\alpha_{j}(nd_{j})}{n^{s}}
   \right)
$$
where $P_{j} = \prod_{i < j} d_{i}$.
\end{lemma}
\noindent This is Lemma 3 of \cite{CGG2} pp.506.

In Lemmas 4 and 5 we consider two
Dirichlet series, $D(s,h/k)$ and $Q(s,\alpha,h/k)$ which  
arise in the analysis of $Q^{*}(s,\alpha,k)$.

\begin{lemma} For $(h,k)=1$ with $k
> 0$ we define
$$
  D \left(s, h/k \right) = \sum_{n =1}^{\infty}
  d(n)n^{-s}e \left( nh/k \right) \ ( \sigma > 1) \ .
$$
Then $D(s,h/k)$ is regular in the entire complex  plane except for
a double pole at $s=1$.  Moreover, it has the same meromorphic part as
$k^{1-2s}\zeta^{2}(s)$.
\end{lemma}
\noindent This is proven in Estermann \cite{E} pp.124-126.

\begin{lemma} Let $(h,k)=1$ and $k = \prod p^{\lambda}
>0$. For $\alpha \in \mathbb{R}$ and $\sigma > 1$ define
\begin{equation}
   \Q \left( s, \alpha , h/k \right)
   = - \sum_{m,n=1}^{\infty} \frac{d(m) \Lambda(n)}{m^{s}
   n^{s-i\alpha}} e \left( \frac{-mnh}{k} \right) \ .
   \label{eq:qsa}
\end{equation}
Then $\Q(s,\alpha,h/k)$ has a meromorphic continuation to the entire complex plane. If $\alpha \ne 0$, $\Q(s,
\alpha,h/k)$ has

\noindent $(i)$ at most a double pole at $s=1$ with same principal
part as
\begin{equation}
   k^{1-2s} \zeta^{2}(s) \left(
   \frac{\zeta^{'}}{\zeta}(s-i \alpha) - \mathcal{G}(s, \alpha,k)
   \right) \ ,
   \label{eq:pol1}
\end{equation}
where
\begin{equation}
   \mathcal{G}(s, \alpha,k) = \sum_{p \mid k} \log p \left(
   \sum_{a=1}^{\lambda-1} p^{a(s-1+i \alpha)}
   + \frac{p^{\lambda(s-1+i \alpha)}}{1-p^{-s+i\alpha}}
   - \frac{1}{p^{s-i \alpha} -1} \right) \ ;
   \label{eq:gsak}
\end{equation}
$(ii)$ a simple pole at $s=1+i \alpha$ with residue
\begin{equation}
   - \frac{1}{k^{i \alpha} \phi(k)} \zeta^{2}(1+i \alpha)
   \rk_{k}(1+i \alpha)
   \label{eq:pol1a}
\end{equation}
where 
\begin{equation}
   \rk_{k}(s) = \prod_{p^{\lambda} \mid  \mid k} (1-p^{-1}+
   \lambda(1-p^{-s})(1-p^{s-1})) \ . 
   \label{eq:rks} 
\end{equation}
Moreover, on GRH, $\mathcal{Q}(s, \alpha,h/k)$ is regular in $\sigma > 1/2$ except for these two poles.
\end{lemma}
\noindent This is Lemma 5 of \cite{CGG} pp.217-218.

In the proof of Lemma 5 of \cite{CGG}, the 
generating function $Q(s,\alpha,h/k)$ is written as a linear combination of 
$(L^{'}/L)(s,\chi)$ where $L(s,\chi)$ is a Dirichlet $L$-function
modulo $k$. These $L$-functions contribute the pole at $s=1+i\alpha$.
Moreover, $Q(s,\alpha,h/k)$ is regular for $\sigma > 1/2$ since
$(L^{'}/L)(s,\chi)$ is regular in this region assuming GRH.

For an arbitrary variable $x$ we define the following generating function for $d_{r}$
\begin{equation}
   T_{r}(x,\lambda) = \sum_{j \ge \lambda} d_{r}(p^{j}) x^{j} \ .
\end{equation}

\begin{lemma} For $r , \lambda \in \mathbb{N}$ and x an indeterminate
we have
\begin{equation}
   (1-x)^{r} T_{r}(x,\lambda) = \lambda d_{r}(p^{\lambda}) \int_{0}^{x} t^{\lambda-1} (1-t)^{r-1} \, dt \ .
   \label{eq:divpoly}
\end{equation}
\end{lemma}
\noindent We define for $\lambda, r \in \mathbb{N}$ the polynomial
\begin{equation}
   H_{\lambda,r}(x) :=  \lambda x^{-\lambda} \int_{0}^{x} t^{\lambda-1}(1-t)^{r-1} \, dt \ .
\end{equation}
Note that $H_{\lambda,r}(x)$ is a degree $r$ polynomial and $H_{\lambda,r}(0)=1$.
Consequently, the lemma may be rewritten as 
\[
   (1-x)^{r} T_{r}(x,\lambda) = d_{r}(p^{\lambda}) x^{\lambda}   H_{\lambda,r}(x)  \ . 
\]
\noindent {\it Proof}. Define the generating
functions
\begin{equation}
\begin{split}
  & A(x,y) := \sum_{\lambda =1}^{\infty} (1-x)^{r} T_{r}(x,\lambda) y^{\lambda} \ , \\
  & B(x,y) := \sum_{\lambda =1}^{\infty}
   \left( \lambda d_{r}(p^{\lambda}) \int_{0}^{x} t^{\lambda-1} (1-t)^{r-1} \, dt \right)
   y^{\lambda} \ .
   \nonumber
\end{split}
\end{equation}
We will show that these generating functions are equal and hence we have~(\ref{eq:divpoly}). Note that
\begin{equation}
\begin{split}
   A(x,y) & = (1-x)^{r} \sum_{j=1}^{\infty} d_{r}(p^{j})x^{j} \sum_{\lambda =1}^{j} y^{\lambda}
   = \frac{y(1-x)^{r}}{y-1} \sum_{j=1}^{\infty} d_{r}(p^{j})x^{j} (y^{j}-1) \\
   & = \frac{y}{y-1}\left( \frac{(1-x)^r}{(1-xy)^{r}}-1  \right)
   \nonumber
\end{split}
\end{equation}
and since $\lambda d_{r}(p^{\lambda}) = rd_{r+1}(p^{\lambda-1})$ for $\lambda \ge 1$
\[
  B(x,y) = r \int_{0}^{x} (1-t)^{r-1} \left( \sum_{\lambda =1}^{\infty} d_{r+1}(p^{\lambda-1}) t^{\lambda-1}
  y^{\lambda} \right) \, dt
   = ry \int_{0}^{x} \frac{(1-t)^{r-1}}{(1-ty)^{r+1}} \, dt \ .
\]
A calculation shows that $A_{x}(x,y)=B_{x}(x,y)=
\frac{ry(1-x)^{r-1}}{(1-xy)^{r+1}}$
and since $A(0,y)=B(0,y)=0$ it follows that $A(x,y)=B(x,y)$.

Our calculations require Perron's formula.
\begin{lemma} Let $F(s) := \sum_{n\ge 1} a_{n} n^{-s}$ be a
Dirichlet series with finite abscissa of absolute convergence
$\sigma_{a}$. Suppose there exists a real number $\alpha \ge 0$
such that
$$
  \sum_{n=1}^{\infty} |a_{n}|n^{-\sigma} \ll (\sigma-\sigma_{a})^{-\alpha}
  \
  (\sigma > \sigma_{a})
$$
and that B is a non-decreasing function such that $|a_{n}| \le
B(n)$ for $n \ge 1$.  Then for $x \ge 2, T \ge 2, \sigma \le
\sigma_{a}, \kappa := \sigma_{a}-\sigma+(\log x)^{-1}$, we have
\begin{equation}
\begin{split}
  \sum_{n \le x} \frac{a_{n}}{n^{s}}
  & = \frac{1}{2 \pi i} \int_{\kappa-iT}^{\kappa+iT} F(s+w)
  \frac{x^{w}}{w} \, dw
  + O \left( \frac{x^{\sigma_{a}-\sigma}(\log x)^{\alpha}}{T}
             + \frac{B(2x)}{x^{\sigma}}
             \left( 1 + x \frac{\log T}{T} \right)
  \right) \ .
\end{split}
\end{equation}
\end{lemma}
\noindent This is Corollary 2.1 p.133 of \cite{Te}.

The following Lemma is another place where GRH is invoked. 
This lemma gives bounds for 
 $Q^{*}(s,\alpha,k)$ in the 
critical strip.  These bounds  are required for estimating 
the left side of the contour in~(\ref{eq:percon}).
In fact, GRH shall be invoked in the form of a Lindel\"{o}f type
bound for Dirichlet $L$-functions.

\begin{lemma} Assume GRH. Let $y = (T/2 \pi)^{\eta}$ where $0 < \eta
< 1/2$, $k \in \mathbb{N}$ with $k \le y$, and $\alpha \in
\mathbb{R}$. Set
\begin{equation}
   \Q^{*}(s, \alpha,k) =
   \sum_{j=1}^{\infty} b(j)j^{-s}e \left( -j/k \right) \
   (\sigma > 1) \ ,
   \label{eq:qstar}
\end{equation}
where
\[
  b(j) = -
   \sum_{{\begin{substack}{hmn=j
         \\ h \le y}\end{substack}}}
  d_{r}(h) P \left( [h]_{y} \right) d(m) \Lambda(n) n^{\i \alpha}
  \ .
\]
Then $\Q^{*}(s, \alpha, k)$ has an analytic continuation to $\sigma
> 1/2$ except possible poles at $s=1$ and $1+i \alpha$.
Furthermore,
\[   \Q^{*}(s, \alpha , k) = O_{\epsilon,P}( y^{\frac{1}{2}}
      T^{\epsilon})
\]
where $s= \sigma +it$,
$\frac{1}{2}+L^{-1} \le \sigma \le 1 + L^{-1}$, $|t| \le T$, $|s-1|> 0.1$,
and $|s-1-i \alpha| > 0.1$.  Note that the constant in the big-$O$ depends on the polynomial $P$. 
\end{lemma}
\noindent {\it Proof}.  If $\chi$ is a character mod $k$, its
Gauss sum is $\tau(\chi)
   = \sum_{h=1}^{k} \chi(h) e \left( h/k
   \right)$ from which it follows that
\begin{equation}
   e \left( - j/k \right)
   = \sum_{d \mid j, d \mid k} \frac{1}{\phi \left( k/d
   \right)} \sum_{ \chi \, (\mathrm{mod} \, \frac{k}{d})}
   \tau(\overline{\chi}) \chi(-j/d) \ .
   \label{eq:ejk}
\end{equation}
By inserting~(\ref{eq:ejk}) in~(\ref{eq:qstar}) we obtain
\begin{equation}
  \Q^{*}(s,\alpha,k) =
  \sum_{d \mid k}\frac{1}{\phi \left( k/d \right) d^{s}}
  \sum_{\chi \, (\mathrm{mod} \, \frac{k}{d})} \tau(\overline{\chi}) \overline{\chi} \left(
  -d \right) B(s,d)    
  \label{eq:qstar2}
\end{equation}
where for $\sigma > 1$, $B(s,d) = \sum_{j=1}^{\infty} b(jd)\chi(jd)j^{-s}$.
We now write 
$P(x)=\sum_{i=0}^{N} c_{i} x^{i}$ and hence we obtain
\begin{equation}
     \Q^{*}(s,\alpha,k) = \sum_{i=0}^{N} \frac{c_{i}}{(\log y)^{i}}
     \Q_{i}^{*}(s,\alpha,k)
     \label{eq:qsum}
\end{equation}
where 
\begin{equation}
      \Q_{i}^{*}(s,\alpha,k) =
  \sum_{d \mid k}\frac{1}{\phi \left( k/d \right) d^{s}}
  \sum_{\chi \, (\mathrm{mod} \, \frac{k}{d})} \tau(\overline{\chi}) \overline{\chi} \left(
  -d \right) \frac{\partial^{i}}{\partial z^{i}} \left. B(s,d;z) \right|_{z=0}    
  \ , 
  \label{eq:qistar}
\end{equation}
\[
   B(s,d;z) = \sum_{j=1}^{\infty}
   b_{z}(dj)\chi(dj) j^{-s}, \ \mathrm{and} \
    b_{z}(j) = 
   \sum_{{\begin{substack}{hmn=j
         \\ h \le y}\end{substack}}}
  d_{r}(h) h^{z} d(m) \Lambda(n) n^{\i \alpha} \ . 
\]
Since $\chi$ is completely multiplicative we note that
\[
   B(s,1;z) = \left( \sum_{h \le y} \frac{\chi(h)d_{r}(h)h^{z}}{h^{s}}
   \right) L(s,\chi)^{2}
   \left( \sum_{n \ge 1} \frac{\chi(n)\Lambda(n)}{n^{s-i \alpha}}
   \right) \ .
\]
An application of Lemma 3 implies
\begin{equation}
  B(s,d;z) = \sum_{f_{1}f_{2}f_{3}f_{4}=d}
           \mathcal{A}_{1}(s,f_{1};z) \mathcal{A}_{2}(s,f_{2},f_{1})
           \mathcal{A}_{2}(s,f_{3},f_{1}f_{2}) \mathcal{A}_{3}(s,f_{4},f_{1}f_{2}f_{3})
  \label{eq:bsd}
\end{equation}
where
\begin{equation}
\begin{split}
   & \mathcal{A}_{1}(s,f;z) = \chi(f) \sum_{h \le y/f}
   \frac{\chi(h)d_{r}(fh)(fh)^z}{h^{s}} \ , \\
  & \mathcal{A}_{2}(s,f,r) = \sum_{(n,r)=1} \frac{\chi(fn)}{n^{s}}
 = \chi(f)L(s,\chi) \prod_{p \mid r} (1- \chi(p)p^{-s}) \ , \\
 & \mathcal{A}_{3}(s,f,r) = - \sum_{(n,r)=1} \chi(fn) \Lambda(fn)
  (fn)^{i \alpha} n^{-s} \ .
\end{split}
\end{equation}
We are aiming to show that uniformly for $|z| \le 0.1L^{-1}$
\begin{equation}
  B(s,d;z) \ll_{\epsilon} \left\{ \begin{array}{cl}
                  y^{\frac{1}{2}} T^{\epsilon}  &  \mbox{if $\chi$ is principal}  \\
                  T^{\epsilon} &  \mbox{otherwise}  \\
                  \end{array}
          \right.
   \label{eq:Bsdz}
\end{equation}
in the region $\sigma \ge 1/2+L^{-1}$, $|t| \le T$, and $|s-1|,|s-1-\i \alpha| > 0.1$.  If~(\ref{eq:Bsdz}) holds then 
we have by applying the Cauchy integral formula with a circle of 
radius $0.1L^{-1}$ that
\begin{equation}
    \frac{\partial^{i}}{\partial z^{i}} \left. B(s,d;z) \right|_{z=0} 
     \ll \left\{ \begin{array}{cl}
                  y^{\frac{1}{2}} T^{\epsilon}  &  \mbox{if $\chi$ is principal}  \\
                  T^{\epsilon} &  \mbox{otherwise}  \\
                  \end{array} 
          \right.   \ . 
     \label{eq:Bder}
\end{equation}
By~(\ref{eq:qistar}) and~(\ref{eq:Bder})
\[
   \Q_{i}^{*}(s,\alpha,k) \ll_{\epsilon} T^{\epsilon}
   \sum_{d \mid k} \frac{1}{\phi(k/d)d^{\frac{1}{2}}}
   \left(
   y^{\frac{1}{2}} |\tau(\chi_{0})| +
   \sum_{\chi \ne \chi_{0} (\mathrm{mod} k/d)} |\tau(\chi)| 
   \right) \ .
\]
Since
\[
  \tau(\chi) \ll \left\{ \begin{array}{cl}
                  (k/d)^{\frac{1}{2}} &  \chi \ne \chi_{0} \, (\mathrm{mod} \ k/d)  \\
                  1 &  \chi = \chi_{0} \, (\mathrm{mod} \ k/d) \\
                  \end{array}
          \right.   
\]
it follows that
\[
  \Q_{i}^{*}(s, \alpha,k) \ll T^{\epsilon} 
  \left(
  (y/k)^{1/2} \sum_{d \mid k} d^{1/2} \phi(d)^{-1} 
  +  k^{1/2}\sum_{d \mid k} d^{-1}
  \right) 
  \ll y^{1/2} T^{\epsilon} 
\]
and hence by~(\ref{eq:qsum})
the desired bound
 $Q^{*}(s,\alpha,k) \ll_{\epsilon,P} yT^{1/2+\epsilon}$ follows.
It now suffices to establish~(\ref{eq:Bsdz}).  If $\chi$ is principal (mod $k/d$) then
\begin{equation}
  \mathcal{A}_{1}(s,f;z) \ll f^{\epsilon} \sum_{n \le y/f}
  n^{-\frac{1}{2}} \ll 
  y^{1/2}  \ .
  \label{eq:a1sfp}
\end{equation}
Now suppose $\chi$ is non-principal.  If $y/f \ll y^{\epsilon}$, we
have trivially that $|\mathcal{A}_{1}(s,f)| \ll y^{\epsilon}$. 
If we suppose $y/f \gg y^{\epsilon}$ then by Perron's formula (Lemma 7)    
\begin{equation}
   \mathcal{A}_{1}(s,f;z) = \frac{\chi(f) f^z}{2 \pi i}
   \int_{\kappa-2iT}^{\kappa+2iT}
   G(s+z+w) \frac{(y/f)^{w}}{w} dw + O
   \left(
   1
   \right)
   \label{eq:a1sf}
\end{equation}
for $\sigma \ge 1/2 + L^{-1}$, $|t| \le T$, 
$\kappa=1-\sigma+2L^{-1}$ where $G(w) = \sum_{n=1}^{\infty}
d_{r}(fn) \chi(n) n^{-w}$.  By multiplicativity we have
\begin{equation}
    G(w) = L(w,\chi)^{r} \prod_{p^{e} \mid \mid f}
  \left(  \frac{  \sum_{a=0}^{\infty} \chi(p^a) d_{r}(p^{e+a})p^{-aw}}{
    \sum_{a=0}^{\infty} \chi(p^{a}) d_{r}(p^{a}) p^{-aw}}  \right) \ . 
    \label{eq:gw}
\end{equation}
By Lemma 6, it follows that
\begin{equation}
   G(w) = d_{r}(f) L(w,\chi)^{r} \prod_{p^\lambda\mid \mid f}
  H_{\lambda,r}(x_p) 
\end{equation}
with $x_p=\chi(p)p^{-s}$.  We have that $|x_p| \le p^{-\sigma}$ and
since $H_{\lambda,r}(0)=1$ it follows that
\[
    \left| \prod_{p^\lambda \mid \mid f}
   H_{\lambda,r}(x_p)   \right|
   \ll \prod_{p \mid \mid f} (1+O(p^{-1/2}) )\ll f^{\epsilon} \ . 
\]
In addition, GRH implies
$|L(w,\chi)| \ll (1+|t|)^{\epsilon} \left( k/d \right)^{\epsilon} $
for $\mathrm{Re}(w) \ge 1/2$ and any $\epsilon > 0$.
We now move the contour in~(\ref{eq:a1sf}) to $\mathrm{Re}(w) = \kappa^{'}$
line where $\kappa^{'}=1/2-\sigma+2L^{-1}$ and we have 
\[
     \mathcal{A}_{1}(s,f;z) = \frac{\chi(f) f^z}{2 \pi i}
   \int_{\kappa^{'}-2iT}^{\kappa^{'}+2iT}
   G(s+z+w) \frac{(y/f)^{w}}{w} dw + O
   \left(
   T^{\epsilon}
   \right)
\] 
Since $0.5 \le \mathrm{Re}(s+z+w)$ and $\mathrm{Re}(w) \le L^{-1}$  it follows that 
\begin{equation}
   \mathcal{A}_{1}(s,f;z)   \ll 
   f^{\epsilon} T^{\epsilon} \left( k/d \right)^{\epsilon}
   \left( y/f \right)^{L^{-1}} 
   \int_{\kappa^{'}-2iT}^{\kappa^{'}+2iT} \frac{|dw|}{|w|}
   \ll T^{\epsilon} \ .
   \label{eq:a1sfnp} 
\end{equation}
For $f$ and $r$ dividing $d$, we have 
\begin{equation}
  \mathcal{A}_{j}(s,f,r) \ll T^{\epsilon}
  \label{eq:ajs}
\end{equation}
for $j=2,3$.  This is proven in \cite{CGG} pp.219-220.
By~(\ref{eq:bsd}) in combination with the bounds~(\ref{eq:a1sfp}),~(\ref{eq:a1sfnp}), and~(\ref{eq:ajs}) we obtain~(\ref{eq:Bsdz}) which finishes
the lemma.

The purpose of the next five lemmas is to provide a variety of formulae 
for mean values of certain multiplicative functions which arise in our asymptotic evaluation of 
$I$~(\ref{eq:Iexp}).  Lemma 9 provides bounds for certain divisor sums.
Lemmas 10, 11, and 13 give asympotic formulae for divisor and other
divisor-like sums.  Lemma 12 provides a formula for simple 
prime number sums.   

\begin{lemma} For $\alpha \in \Bbb R$ and $j \in \Bbb Z_{\ge 0}$ we have
\begin{equation}
   \mathcal{G}^{(j)}(1,\alpha,k) = \sum_{p \mid k} p^{i \alpha} (\log
   p)^{j+1} + O(C_{j}(k))
  \label{eq:gj1}
\end{equation}
where $\mathcal{G}(s,\alpha,k)$ is defined by~(\ref{eq:gsak}) and
\begin{equation}
   C_{j}(k) = \sum_{p \mid k} \frac{\log^{j} p}{p}
            + \sum_{p^{a} \mid \mid k, \ a \ge 2} a \log^{j} p \ .
  \label{eq:Cjk}
\end{equation}
Moreover, we have
\begin{equation}
   \sum_{h,k \le x} \frac{d_{r}(h)d_{r}(k)(h,k)}{hk}
   C_{j} \left( \frac{k}{(h,k)} \right) \ll (\log x)^{r^2+r}   \ .
   \label{eq:avcj}
\end{equation}
\end{lemma}
\noindent {\it Proof}.  We remark that~(\ref{eq:gj1}) is proven in \cite{CGG} pp.222-223.  The sum
in~(\ref{eq:avcj}) is bounded by
\begin{equation}
\begin{split}
  & \sum_{h,k \le x} \frac{d_{r}(h)d_{r}(k)}{hk} (C_{j}(k)+1)
    \sum_{{\begin{substack}{a \mid h
         \\ a \mid k}\end{substack}}}
   \phi(a) \\
  &
  \le \sum_{a \le x} \frac{d_{r}(a)^2\phi(a)}{a^{2}}
  \sum_{h,k \le \frac{x}{a}} \frac{d_{r}(h)d_{r}(k)(C_{j}(ak)+1)}{hk}  \\
  & \le (\log x)^{2r} \sum_{a \le x} \frac{d_{r}(a)^{2}(C_{j}(a)+1)}{a} +
  (\log x)^{r} \sum_{a \le x} \frac{d_{r}(a)^{2}}{a} \sum_{k \le \frac{x}{a}} \frac{d_{r}(k)C_{j}(k)}{k}  
   \nonumber
\end{split}
\end{equation}
Observe that
\begin{equation}
\begin{split}
    \sum_{a \le y} \frac{d_{r}(a)^{2}C_{j}(a)}{a}
    & = \sum_{p \le y} \frac{\log^{j} p}{p}
    \sum_{u \le \frac{y}{p}} \frac{d_{r}(up)^{2}}{up}
    + \sum_{p^{a} \le y, a \ge 2} a(\log p)^{j}
    \sum_{u \le \frac{y}{p^{a}}} \frac{d_{r}(up^{a})^{2}}{up^{a}} \\
    & \ll (\log x)^{r^2} \left( \sum_{p} \frac{(\log p)^{j}}{p^2}  \right) \ll (\log x)^{r^2} 
    \nonumber
\end{split}
\end{equation}
where we have applied~(\ref{eq:ram}). 
A similar argument establishes that $\sum_{k \le  x} d_{r}(k)C_{j}(k)k^{-1} \ll (\log x)^{r^2+r}$ . Putting together the results establishes the lemma.

We now introduce the arithmetic function $\sigma_{r}(m,s)$ where
$r \in \mathbb{N}$ and $s \in \mathbb{C}$. It is defined by
\begin{equation}
  \sigma_{r}(m,s) := \left( \sum_{n=1}^{\infty} \frac{d_{r}(mn)}{n^{s}} 
  \right)   \zeta(s)^{-r}  =
  \prod_{p^{\lambda} \mid \mid m} (1-p^{-s})^{r} p^{\lambda s}
   \sum_{j \ge \lambda}^{\infty} \frac{d_{r}(p^{j})}{p^{js}} \ . 
  \label{eq:sigms0}
\end{equation}
The second equation is obtained by mutiplicativity.  By Lemma 6, it 
follows that 
\[
   \sigma_{r}(m,s) = \prod_{p^{\lambda} \mid \mid m}
   d_{r}(p^{\lambda}) H_{\lambda,r}(p^{-s}) \ . 
\]
The value $s=1$ will have a special importance so we set
$\sigma_{r}(m) := \sigma_{r}(m,1)$.
In the following calculations we shall often employ the bound 
\begin{equation}
  |\sigma_{r}(m,s)| \ll d_{r}(m) j_{\tau}(m)  \ \mathrm{for} \ 
  \mathrm{Re}(s) \ge \tau > 0
  \label{eq:sigrbd}
\end{equation}
The function $\sigma_{r}$  is a correction factor that arises  
due to the fact $d_{r}$ is not completely mutiplicative.  More precisely, 
we notice in all cases of the following lemma that
\[
   \sum_{h \le t} d_{r}(mh) f(h) \sim \sigma_{r}(m) \sum_{h \le t} d_{r}(h) f(h)
\]
where $f$ is a smooth function.

\begin{lemma} Suppose $r,n \in \mathbb{N}$, $1 \le x,n \le
\frac{T}{2 \pi}$, and $F \in C^{1}([0,1])$.  There exists an
absolute constant $\tau_{0} = \tau_{0}(r)$ such  that
\begin{equation}
   \sum_{h \le x} \frac{d_{r}(nh)}{h} F \left( [h]_{x}
   \right)
   = \frac{\sigma_{r}(n)(\log x)^{r}}{(r-1)!}
   \int_{0}^{1} \theta^{r-1} F(\theta) \, d \theta
  + O(d_{r}(n)j_{\tau_0}(n)) 
\end{equation}
where $j_{\tau_{0}}(n)$ is defined by~(\ref{eq:jtau}).  In order to abbreviate
notation we define 
\begin{equation}
  \epsilon(n) = d_{r}(n) j_{\tau_{0}}(n)   \ . 
  \label{eq:eumt}
\end{equation}
Suppose $m,u,v \in \mathbb{N}$, $1 \le y,m \le
\frac{T}{2 \pi}$, $p$ a prime with $p \le \frac{T}{2 \pi}$,
and $P \in C^{1}([0,1])$.
 We now deduce the following formulae:  \\
$(i)$
\begin{multline}
  \sum_{h \le \frac{y}{m}} \frac{d_{r}(mh)}{h} (\log h)^{u}
  P([mh]_{y})
  \sim \frac{\sigma_{r}(m)}{(r-1)!} \log \left( \frac{y}{m} \right)^{r+u}
  \int_{0}^{1} F_{1}(\theta,m) \, d \theta  \ ,
  \label{eq:div1}
\end{multline}
$(ii)$
\begin{multline}
  \sum_{h \le \frac{y}{mp}} \frac{d_{r}(mph)P([mph]_{y})(\log
             ph)^{v}}{h}
   \sim \frac{\sigma_{r}(pm)}{(r-1)!}
   \log \left( \frac{y}{pm} \right)^{r}
   \int_{0}^{1}
   F_{2}(\theta,pm) \, d \theta   \ ,
   \label{eq:div2}
\end{multline}
$(iii)$
\begin{multline}
   \sum_{h \le \frac{y}{m}} \frac{d_{r}(mh)P([mh]_{y})}{h}
  \log \left( \frac{T}{2 \pi h} \right)^{u}
  \sim \frac{\sigma_{r}(m)(\log y)^{u+r}}{(r-1)!}
  \int_{0}^{1-[m]_{y}} F_{3}(\theta,m)  \, d\theta \ ,
  \label{eq:div3}
\end{multline}
where each formula is valid up to an error term 
$\epsilon(m)=d_{r}(m) j_{\tau_0}(m)$ and
\begin{equation}
\begin{split}
  & F_{1}(\theta,m) = \theta^{r+u-1} P \left( [m]_{y}+(1-[m]_{y})\theta  \right)  \ , \\
  & F_{2}(\theta,pm) = \theta^{r-1}
   \left(
   \log p + \theta \log \frac{y}{pm}
   \right)^{v}
   P\left(
   [pm]_{y} + (1-[pm]_{y}) \theta
   \right) \ , \\
   & F_{3}(\theta,m) = \theta^{r-1}(\eta^{-1}-\theta)^{u}
  P \left( [m]_{y} + \theta \right) \ .
  \label{eq:Fi}
\end{split}
\end{equation}
\end{lemma}

\noindent {\it Proof}.  It was established in Lemmas 4 and 5 of \cite{CG} that
\begin{equation}
   \sum_{h \le t} \frac{d_{r}(nh)}{h} = \frac{\sigma_{r}(n)(\log t)^{r}}{r!} +
   O(d_{r}(n) j_{\tau_{0}}(n))
   \label{eq:sel}
\end{equation}
for some $\tau_{0}=\tau_{0}(r) > 0$. We abbreviate~(\ref{eq:sel}) to $T(t) = M(t) + O(\epsilon(n))$. If $g \in
\mathcal{C}^{1}([0,1])$ we deduce
\begin{multline}
   \sum_{h \le x} \frac{d_{r}(nh)}{h} g \left( [h]_{x}
   \right)
   = \int_{1}^{x} M^{'}(t) g\left(  [t]_{x} \right) \, dt \\
    + O \left(\epsilon(n) \left(|g(0)|+|g(1)| + \frac{1}{\log x}
   \int_{1}^{x} |g^{'} \left( [t]_{x} \right)| \, \frac{dt}{t} \right) \right) \ .
   \nonumber
\end{multline}
The error term is $\ll \epsilon(n)$ and the principal term is
\begin{equation}
   \frac{\sigma_{r}(n)}{(r-1)!} \int_{1}^{x} (\log t)^{r-1}  g\left(
   [t]_{x} \right) \, dt
   =  \frac{\sigma_{r}(n)(\log x)^{r-1}}{(r-1)!}
   \int_{0}^{1} \theta^{r-1} g(\theta) \, d \theta
   \nonumber
\end{equation}
by the variable change $\theta = [t]_{x}$. Formulae $(i)$-$(iii)$ of this lemma correspond to the following
choices of parameters $(n,g(\theta),x)$:
\[
  \left( m \ , \theta^{u}P([m]_{y}+\theta), \frac{y}{m}   \right) \ ,
  \left( pm \ , \left( [p]_{x} + \theta
  \right)^{u}P([pm]_{y} + \theta),
  \frac{y}{pm} \right) \ ,
\]
\[
  \left( m \ , \left( \frac{\log \frac{T}{2 \pi}}{\log x}
  - \theta \right)^{u}P([m]_{y}+\theta)  ,  \frac{y}{m} \right) \ .
\]
Note that the error term in $(ii)$ is $\epsilon(pm) \ll \epsilon(m)$.  Furthermore, part $(iii)$ requires the
variable change $\theta \to [x]_{y} \theta$.

In the following lemma we consider averages of the expression $\sigma_{r}(\cdot)^{2}$.  It is in this lemma that the constant 
$a_{r+1}$~(\ref{eq:ar}) of Theorem 2 appears.  It naturally arises   
upon considering the Dirichlet series 
$\sum_{n \ge 1} \phi(n) \sigma_{r}(n)^2 n^{-s}$.

\begin{lemma}
Let $r \in \N$ and $g \in C^{1}([0,1])$.  \\
 $(i)$ For $p \le y$ prime we have 
\begin{equation}
\begin{split}
  \sum_{m \le \frac{y}{p}} \frac{\phi(m) \sigma_{r}(m) \sigma_{r}(pm) }{m^{2}}
  g \left( [m]_{y} \right)
 & =  \frac{\sigma_{r}(p)a_{r+1}(\log y)^{r^2}}{(r^2-1)!}
 \int_{0}^{1-[p]_{y}} \delta^{r^2-1} g(\delta) \, d \delta \\
 & + O \left( (\log y)^{r^2}(p^{-1} + (\log y)^{-1} ) \right)  \ .
 \label{eq:sig1}
\end{split}
\end{equation}

\noindent $(ii)$ For $0 \le \theta < 1 $ we have
\begin{equation}
   \sum_{m \le y^{1-\theta}} \frac{\phi(m) \sigma_{r}(m)^{2}}{m^{2}}
   g \left( [m]_{y} \right)
   = \frac{a_{r+1}(\log y)^{r^2}}{(r^2-1)!}
   \int_{0}^{1- \theta} \delta^{r^2-1} g(\delta) \, d \delta (1+O((\log y)^{-1})) \ .
   \label{eq:sig2}
\end{equation}
\end{lemma}
\noindent {\it Proof}.  We only prove $(i)$ since $(ii)$ is similar. We begin by noting that
\begin{equation}
\begin{split}
   \sum_{m \le t} \frac{\phi(m) \sigma_{r}(m) \sigma_{r}(pm)}{m^2}
   & = \sigma_{r}(p) \sum_{m \le t} \frac{\phi(m) \sigma_{r}(m)^{2}}{m^2} \\
   & +  \sum_{{\begin{substack}{m \le t
         \\ p \mid m}\end{substack}}}
   \frac{\phi(m)\sigma_{r}(m)(\sigma_{r}(p) \sigma_{r}(m)-\sigma_{r}(pm))}{m^2} \ .
    \nonumber
\end{split}
\end{equation}
Since, $\sigma_{r}(m) \ll d_{r}(m)j_{1}(m)$, $d_{r}(uv) \le d_{r}(u)d_{r}(v)$, and $\phi(up) \le \phi(u) p$, it
follows that the second term is
\[
 \ll \frac{d_{r}(p)}{p} \sum_{n \le \frac{x}{p}} \frac{d_{r}(n)^2 j_{1}(n)}{n} \ll p^{-1} (\log x)^{r^2} \ .
\]
By equations (36)-(38) of \cite{CG} in conjunction with Theorem 2 of \cite{Se} we deduce
\[
     \sum_{m \le t}
  \frac{\phi(m) \sigma_{r}(m)^{2}}{m^2}
  = \frac{a_{r+1}(\log t)^{r^2}}{(r^2-1)!} (1+ O((\log t)^{-1}))
\]
and hence we arrive at
\[
     \sum_{m \le t} \frac{\phi(m) \sigma_{r}(m) \sigma_{r}(pm)}{m^2}  =
  \frac{ra_{r+1}(\log t)^{r^{2}}}{(r^{2}-1)!} +
    O \left( (\log t)^{r^{2}}p^{-1} + \log^{r^{2}-1} t \right) \ .
\]
We abbreviate this equation to $T(t) = M(t) + O(E(t))$.  The sum in $(i)$ may be expressed as the Stieltjes
integral
\begin{equation}
\begin{split}
  \int_{1^{-}}^{\frac{y}{p}} g \left( [t]_{y} \right) dT(t)
  = \int_{1^{-}}^{\frac{y}{p}} g \left( [t]_{y} \right) dM(t)
  +  \left. g \left( [t]_{y} \right) E(t) \right|_{1^{-}}^{\frac{y}{p}}
  - \int_{1^{-}}^{\frac{y}{p}} g^{'} \left( [t]_{y} \right) E(t) dt \ .
  \label{eq:sigmH}
\end{split}
\end{equation}
The integral equals
\[
\frac{ra_{r+1}}{(r^2-1)!} \int_{1}^{\frac{y}{p}} (\log t)^{r^2-1}
  g \left( [t]_{y} \right) \, dt =
  \frac{ra_{r+1}}{(r^2-1)!} \int_{0}^{1-[p]_{y}}
  \delta^{r^2-1} g(\delta) \, d\delta \ .
\]
Moreover, it is clear that the error term in~(\ref{eq:sigmH}) is $O \left( (\log t)^{r^{2}} p^{-1} + (\log
t)^{r^2-1} \right)$.

In the main calculation of this article we compute certain simple sums over primes.  The following lemma provides
the required result

\begin{lemma} Suppose $w \ge 1$, $0 \le \theta < 1$, and $g \in C^{1}([0,1])$
then
\begin{equation}
\begin{split}
  \sum_{p \le y^{1-\theta}} \frac{(\log p)^{w}}{p^{1-i \alpha}}
  g\left( [p]_{y} \right)
  & = \sum_{j=0}^{\infty} \frac{(i \alpha)^{j}}{j!} (\log
   y)^{j+w}
   \int_{0}^{1-\theta} \beta^{j+w-1} g(\beta) \, d\beta \\
   & + O( (\log y)^{w-1}) \ .
\end{split}
\end{equation}
\end{lemma}
\noindent {\it Proof}. By Stieltjes integration the sum in question is
$$
  \int_{1}^{y^{1-\theta}} t^{i \alpha} (\log t)^{w-1}
  g\left( [t]_{y} \right) \frac{d\theta(t)}{t}
$$
where $\theta(t) = \sum_{p \le t} \log p = t+\epsilon(t)$ and $\epsilon(t) \ll t \exp(-c \sqrt{\log t})$. Note
that the main term is
$$
  \int_{1}^{y^{1-\theta}} t^{i \alpha} (\log t)^{w-1}
  g\left( [t]_{y} \right) \frac{dt}{t}
  = \sum_{j=0}^{\infty} \frac{(i \alpha)^{j}}{j!}
  \int_{1}^{y^{1-\theta}} (\log t)^{j+w-1}
  g\left( [t]_{y} \right) \frac{dt}{t} \ .
$$
By the variable change $\beta=[t]_{y}$ we obtain the required expression for the principal part.  Put $h(t) =
t^{i\alpha} (\log t)^{w-1} g \left( [t]_{y} \right) t^{-1}$ and note $h(t) \ll (\log t)^{w-1}t^{-1}$ and $h^{'}(t)
\ll (\log t)^{w-1}t^{-2}$ for $t \le y$. By the above bound for $\epsilon(t)$
\[
   \int_{1}^{y^{1-\theta}} h(t) d \epsilon(t)
   \ll h(y^{1-\theta}) \epsilon(y^{1-\theta})
   + \int_{1}^{y^{1-\theta}} h^{'}(t) \epsilon(t) dt \ll (\log y)^{w-1}  \ . 
\]

We now define
$f(k) =\rk_{k}(1+i \alpha)/\phi(k)$ where
$\rk_{k}(s)$ is defined by~(\ref{eq:rks}). 
In the following lemmas we shall study the Dirichlet series 
\begin{equation}
  Z(s,\alpha) = \sum_{k \ge 1} d_{r}(mk)f(nk) k^{-s} 
  =  \sum_{k \ge 1} \frac{d_{r}(mk) \rk_{nk}(1+i \alpha)}{\phi(nk) k^{s}}  \ . 
  \label{eq:zsalpha}
\end{equation}
Since $f$ is a multiplicative function, it is determined by its value
at the prime powers.  Consequently, we could equivalently define $f$ by the rule
\begin{equation}
   f(p^{a}) := (1+ak_{p})p^{-a}
   \label{eq:fp}
\end{equation}
where
\begin{equation}
   k_{p} := k_{p}(\alpha) = (1-p^{i \alpha})(1-p^{-1-i \alpha})/(1-p^{-1}) 
   \label{eq:kp}
\end{equation}
which we obtain from~(\ref{eq:rks}).  Moreover, note that $k_{p}(0) = 0$.
\begin{lemma} Put $l = \log x$ and suppose $|\alpha|
\ll (\log x)^{-1}$. For $1 \le m \le T$, $n$ squarefree and $n \mid m$
we have 
\[
   \sum_{k \le x} d_{r}(mk)f(nk) =
   \frac{\sigma_{r}(m)}{n} l^{r}
   \sum_{j=0}^{r} \binom{r}{j} \frac{1}{(r+j)!} (-i \alpha \ l)^{j} +
   O \left( \frac{d_{r}(m)j_{\tau_{0}}(m)l^{r-1}}{n^{1-\epsilon}} \right)
\]
where $\tau_{0}=1/3$ is valid and $j_{\tau_{0}}(m)$ 
is defined by~(\ref{eq:jtau}).
\end{lemma}
\noindent {\it Proof}.  This lemma will follow from an application of Perron's formula.
However, we must begin by analyzing the Dirichlet series $Z(s,\alpha)$.  We put $m = \prod_{p} p^{\lambda}=uv$ with $u = \prod_{p \mid n} p^{\lambda}$
and
hence by multiplicativity 
\begin{equation}
  Z(s,\alpha) =
         \left(   \prod_{p^{\lambda} \mid \mid u}   
         \frac{\alpha_{p}(s,\alpha)}{h_{p}(s,\alpha)}  \right)
         \prod_{p^{\lambda}
          \mid \mid v} \left( \frac{\beta_{p}(s,\alpha)}{h_{p}(s,\alpha)} \right)
         \left( \prod_{p} h_{p}(s,\alpha) \right)
\end{equation}
where
\begin{equation}
  \alpha_{p} = \alpha_{p}(s,\alpha) = \sum_{a \ge 0}
  d_{r}(p^{a+\lambda})f(p^{a+1})p^{-as} \ ,
  \label{eq:alps}
\end{equation}
\begin{equation}
  \beta_{p} = \beta_{p}(s,\alpha) = \sum_{a \ge 0}
  d_{r}(p^{a+\lambda})f(p^{a})p^{-as} \ ,
  \label{eq:beps}
\end{equation}
\begin{equation}
  h_{p} = h_{p}(s,\alpha) = \sum_{a \ge 0}
  d_{r}(p^{a})f(p^{a})p^{-as} \ . \
  \label{eq:hps}
\end{equation}
In the above product we label
\begin{equation}
  Z_{11}(s,\alpha) = \prod_{p^{\lambda} \mid \mid u} \frac{\alpha_{p}(s,\alpha)}{h_{p}(s,\alpha)} \ , \
  Z_{12}(s,\alpha) = \prod_{p^{\lambda} \mid \mid v} \frac{\beta_{p}(s,\alpha)}{h_{p}(s,\alpha)} \ ,
  \label{eq:z11z12}
\end{equation}
and we set $Z_{1}(s,\alpha)= Z_{11}(s,\alpha)Z_{12}(s,\alpha)$.  Next we
remark that the last product factors as  
\begin{equation}
  \prod_{p} h_{p}(s,\alpha) = \frac{\zeta^{2r}(1+s)}{\zeta^{r}(1+s-i\alpha)} Z_{3}(s,\alpha)
  := Z_{2}(s,\alpha) Z_{3}(s,\alpha)
  \label{eq:z2z3}
\end{equation}
with $Z_{3}(s,\alpha)$ holomorphic in $\mathrm{Re}(s) > -1/2$. 
This shall follow from the expressions we derive for $\alpha_p,\beta_p,$
and $h_p$ in the next section.  Thus we have the factorization
\begin{equation}
  Z(s,\alpha) = Z_{1}(s,\alpha)Z_{2}(s,\alpha)Z_{3}(s,\alpha) \ .
  \label{eq:zsa}
\end{equation}
By Perron's formula we have 
\begin{equation}
\begin{split}
  \sum_{k \le x} d_{r}(mk)f(nk)
  & = \frac{1}{2 \pi i} \int_{c-iU}^{c+iU} Z(s,\alpha) \frac{x^{s}}{s} \,
  ds + O  \left( \frac{d_{r}(m)}{n^{1-\epsilon}}
  \left(\frac{(\log x)^{2r}}{U} + 1   \right)      \right) 
  \label{eq:perr}
\end{split}
\end{equation}
where $c=(\log x)^{-1}$.  
Let $\Gamma(U)$ denote the contour 
consisting of $s \in \mathbb{C}$ such that 
\[
    \mathrm{Re}(s)=- \frac{\beta}{\log ( |\mathrm{Im}(s)| + 2)}
\]
where $\beta$ is a sufficiently small fixed positive number and $|\mathrm{Im}(s)| \le U$.  Our strategy will be to deform the contour in~(\ref{eq:perr}) to $\Gamma(U)$, thus picking up the pole at $s=0$ which shall account 
for the main term in the lemma.  However, we must also bound the 
contribution coming from $\Gamma(U)$ and the horizontal parts of 
the contour.   In the following section, we shall establish 
\begin{equation}
    |Z_{1}(s,\alpha)| \ll \frac{d_{r}(m) j_{\tau_{0}}(m)}{n^{1-\epsilon}}
    \label{eq:z1s}
\end{equation}
in the cases $\mathrm{Re}(s) \ge -1/2$, $|\alpha| \le cL^{-1}$ and also 
$\mathrm{Re}(s) \ge -\epsilon$, $|\alpha| \le \epsilon$.  
Moreover, we have $|Z_{3}(s,\alpha)| \ll 1$ in $\mathrm{Re}(s) \ge -1/4$ by the absolute convergence of its series.  Furthermore, it is known that 
\[
      \zeta(1+s)-\frac{1}{s} = O ( \log ( |\mathrm{Im}(s)| + 2))   
      \ \mathrm{and} \
      \frac{1}{\zeta(1+s)}  = O ( \log ( |\mathrm{Im}(s)| + 2))
\]
on $\Gamma(U) $ and to the right of $\Gamma(U)$. 
By~(\ref{eq:zsa}) and our previous estimates, we have 
on $\Gamma(U)$ the bound  
\begin{equation}
    |Z(s,\alpha)| \ll  \log ( |\mathrm{Im}(s)| + 2)^{3r} \frac{d_{r}(m)j_{\tau_{0}}(m)}{n^{1-\epsilon}} 
    \  .
    \label{eq:zsal}
\end{equation}   
We now deform the above contour to $\Gamma(U)$ picking up 
the residue at $s=0$. It follows that 
\begin{equation}
\begin{split}
   & \frac{1}{2 \pi i} \int_{\Gamma(U)} 
   Z(s,\alpha) \frac{x^{s}}{s}  \, ds
   \ll  \frac{d_{r}(m)j_{\tau_{0}}(m)}{n^{1-\epsilon}} 
   \int_{0}^{U} x^{-\frac{\beta}{\log(|t|+2)} }
   (\log (t+2))^{3r} \frac{dt}{|t|+1} \\
   &  \ll \frac{d_{r}(m)j_{\tau_{0}}(m)}{n^{1-\epsilon}} (\log U)^{3r+1}
   \exp \left(-\frac{\beta \log x}{\log(U+2)} \right) \\
   & \ll \frac{d_{r}(m)j_{\tau_{0}}(m)}{n^{1-\epsilon}}
   \exp(-\beta_1 \sqrt{\log x})
   \label{eq:deform}
\end{split}
\end{equation}
by the choice $U = \exp(\beta_2 \sqrt{\log x})$ for a suitable $\beta_2$. 
Similarly, we can show that the horizontal edges connecting $\Gamma(U)$
to $[c-iU,c+iU]$ contribute an amount $d_{r}(m)j_{\tau_0}(m)n^{\epsilon-1}
U^{\epsilon-1}$.  Collecting estimates we conclude 
\begin{equation}
\begin{split}
  \sum_{k \le x} d_{r}(mk)f(nk)
  & =  {\begin{substack} { res \\ s=0} \end{substack}} \left(Z(s,\alpha)x^{s}s^{-1} \right) +
  O ( d_{r}(m)j_{\tau_{0}}(m)n^{\epsilon-1}  )   \ .
  \label{eq:main}
\end{split}
\end{equation}
 
In the next two subsections we establish the bound~(\ref{eq:z1s})
and in the final subsection we will compute the residue in~(\ref{eq:main}). 
\subsection{Computing the local factors $h_{p}$, $\alpha_{p}$, and $\beta_{p}$}
We simplify notation by putting $u=p^{-s-1}$ and $s = \sigma +i t$.  By~(\ref{eq:fp}) 
and~(\ref{eq:hps}) we have
\begin{equation}
   h_{p} = \sum_{a=0}^{\infty} d_{r}(p^{a}) u^{a} + k_{p} \sum_{a=0}^{\infty} a d_{r}(p^{a}) u^{a}
   =   (1-u)^{-r-1}(1 + (rk_{p}-1)u) \ .
   \label{eq:hps}
\end{equation}
Note that we have use $a d_{r}(p^{a}) = r d_{r+1}(p^{a-1})$ for $a \ge 1$.
By~(\ref{eq:kp}), $k_{p}= 1- p^{i \alpha} + O (p^{-1+\epsilon})$
and it follows that
\begin{equation}
       h_{p} =(1-p^{-s-1})^{-r-1}
         \left( 1 + \frac{r-1}{p^{s+1}} - \frac{r}{p^{s+1-i \alpha}} + O(p^{-2-\sigma + \epsilon}) \right) \ .
         \label{eq:hps2} \
\end{equation}
Equation~(\ref{eq:z2z3}) now follows from~(\ref{eq:hps2}). As before we have for $\lambda \ge 1$
\[
  \beta_{p} = \sum_{a=0}^{\infty} d_{r}(p^{a+\lambda}) u^{a} + k_{p} \sum_{a=0}^{\infty} ad_{r}(p^{a+\lambda}) u^{a}
  := \beta + k_{p} \tilde{\beta}  \ .
\]
Note that by Lemma 6, $\beta =d_{r}(p^{\lambda})(1-u)^{-r} H_{\lambda,r}(u)$.
and hence it follows that 
\[
  \beta = d_{r}(p^{\lambda}) (1-u)^{-r-1}(1+O_{r}(p^{-1-\sigma})) \ .
\]
Similarly, we note that $\tilde{\beta} = u \frac{d}{du}(\beta(u))$ from which it follows 
that
\[
    \tilde{\beta} =  d_{r}(p^{\lambda})u(1-u)^{-r-1}
     ( (1-u)\frac{d}{du} H_{\lambda,r}(u)-rH_{\lambda,r}(u))
     \ll d_{r}(p^{\lambda})(1-u)^{-r-1} O(|u|) \ . 
\]
We conclude that
\begin{equation}
   \beta_{p}    = d_{r}(p^{\lambda})(1-u)^{-r-1}
   \left(1  + O \left( |k_{p}| p^{-1-\sigma} \right) \right) \ .
   \label{eq:beps2}
\end{equation}
Likewise, we have
\begin{equation}
\begin{split}
   \alpha_{p}   & =  \frac{1}{p} \left(  \sum_{a=0}^{\infty} d_{r}(p^{a+\lambda}) u^{a} +
    k_{p} \left( u \sum_{a=0}^{\infty} d_{r}(p^{a+\lambda})u^{a}   \right)^{'}  \right) 
   =  \frac{1}{p} \left( \beta(1+k_{p}) + k_{p}  \tilde{\beta} \right) \ 
   \nonumber
\end{split}
\end{equation}
and it follows from our previous estimates that
\begin{equation}
   \alpha_{p} = d_{r}(p^{\lambda}) p^{-1}
   (1-u)^{-r-1} O( 
  ( |k_{p}| + 1)) \ . 
  \label{eq:alps2}
\end{equation}

\subsection{Establishing~(\ref{eq:z1s})}

With our estimates for $\alpha_p,\beta_p$, and $h_p$ in hand, we 
are ready to estimate $Z_{1i}(s,\alpha)$.  We have 
by~(\ref{eq:z11z12}),~(\ref{eq:hps}), and~(\ref{eq:alps2})
\begin{equation}
    |Z_{11}(s,\alpha)| \le 
    \prod_{p^{\lambda} \mid \mid u} \frac{|\alpha_p(s,\alpha)|}{|h_p(s,\alpha)|}
    \le 
     \prod_{p^{\lambda} \mid \mid u} 
     \frac{d_{r}(p^{\lambda}) (|k_p|+1)}{p}
     |1+(rk_{p}-1)p^{-s-1}|^{-1} \ . 
     \label{eq:z11b}
\end{equation}
In addition, by~(\ref{eq:z11z12}),~(\ref{eq:hps}), and~(\ref{eq:beps2}) 
it follows that 
\begin{equation}
    |Z_{12}(s,\alpha)| \le 
    \prod_{p^{\lambda} \mid \mid v} \frac{|\beta_p(s,\alpha)|}{|h_p(s,\alpha)|}
    \le 
     \prod_{p^{\lambda} \mid \mid v} 
     d_{r}(p^{\lambda}) (1+ O(|k_p|p^{-1-\sigma}))
     |1+(rk_{p}-1)p^{-s-1}|^{-1} \ . 
     \label{eq:z12b}
\end{equation}
In order to finish bounding these terms, we require a bound for $k_p$.  
We shall provide a bound for $k_p$ and hence $Z_{1i}(s,\alpha)$ in
each of the cases
$0 < |\alpha| \le cL^{-1}$ and $0 < |\alpha| \le \epsilon$.    

\noindent {\bf Case 1}: $0 < |\alpha| \le cL^{-1}$ and $\mathrm{Re}(s) \ge -1/2$. \\
By the definition~(\ref{eq:kp}) it follows that
\begin{equation}
    |k_p| \ll_{c} |1-p^{i  \alpha}| \ll_{c} \min \left(1,\frac{\log p}{L} \right)
    \label{eq:kpb}
\end{equation}
since we have the bounds $|p^{i \alpha}| \le \exp( |\alpha| \log p)$ and 
 $|1-p^{i \alpha}| \le (|\alpha| \log p) e^{|\alpha| \log p}$. 
 Let $c_{1},c_{2}, \ldots $ be effectively computable constants depending on $c$ and $r$.  We have
$|(rk_{p}-1)p^{-s-1}| \ll p^{-\frac{1}{2}} < 0.5$
if $p \ge c_{1}$. If $p \le c_{1}$ then we may choose $T$ 
sufficiently large such that~(\ref{eq:kpb}) yields $|k_{p}| \le 1/20r$.  Thus $|(rk_{p}-1)p^{-s-1}| \le 1.1p^{-\frac{1}{2}} < 0.8$
for all primes $p < c_{1}$ as long as $T$ is sufficiently large.  By~(\ref{eq:z11b})
and our aforementioned bounds we obtain,
\begin{equation}
  |Z_{11}(s,\alpha)| \le \prod_{p^{\lambda} \mid \mid u} \frac{c_{2}d_{r}(p^{\lambda})}{p} 
  \le \frac{d_{r}(u) c_{2}^{\nu(n)}}{n}
\end{equation}
where $\nu(n)$ is the number of prime factors of $n$ and
\begin{equation}
  Z_{12}(s,\alpha) = \prod_{p^{\lambda} \mid \mid v} d_{r}(p^{\lambda})
   \left( \frac{1 + O(p^{-\frac{1}{2} + \epsilon})}{1 + O(p^{-\frac{1}{2} + \epsilon})} \right)
   = d_{r}(v) \prod_{p \mid v} (1 + O(p^{-1/2+\epsilon})) \ .
\end{equation}
Since $c_{2}^{\nu(n)} \ll n^{\epsilon}$ and 
$Z_1(s,\alpha)=Z_{11}(s,\alpha)Z_{12}(s,\alpha)$ we deduce that 
$Z_{1}(s,\alpha) \ll  d_{r}(m)j_{1/3}(v) n^{\epsilon-1}$
in the range $\mathrm{Re}(s) \ge -1/2$ and $|\alpha| \le cL^{-1}$.

\noindent {\bf Case 2}: $0 < |\alpha| \le \epsilon$ and $\mathrm{Re}(s) \ge -\epsilon$. \\
In this case, it follows from~(\ref{eq:kp}) that 
\begin{equation}
  |k_{p}| \le 4 |1-p^{i \alpha}| \le 
  \min \left(8p^{\epsilon},4\epsilon (\log p) p^{\epsilon} \right) 
  \label{eq:kpb2} 
\end{equation}
by employing again the bounds $|p^{i \alpha}| \le \exp( |\alpha| \log p)$ and 
 $|1-p^{i \alpha}| \le (|\alpha| \log p) e^{|\alpha| \log p}$. 
The first bound in~(\ref{eq:kpb2}) implies that 
$
   |(rk_{p}-1)p^{-s-1}| \le (8r+1)p^{-1+2\epsilon} < 0.5
$
if $p$ is sufficiently large, say $p > c_{3}$.  If $p \le  c_{3}$ then
$
    |(rk_{p}-1)p^{-s-1}|  \le \frac{4r \epsilon (\log p)}{p^{1-\epsilon}}
    + \frac{1}{2^{1-\epsilon}} \le  0.51$ for $\epsilon$ sufficiently small.  Thus
\[
   Z_{11}(s,\alpha) = \prod_{p^{\lambda} \mid \mid u,p \le c_{3}}
    \left(\frac{c_{4}d_{r}(p^{\lambda})}{p^{1-\epsilon}} \right)
   \prod_{p^{\lambda} \mid \mid u, p > c_{3}} \frac{d_{r}(p^{\lambda})}{p^{1-\epsilon}}(1 + O(p^{-1+\epsilon}))
   \ll \frac{d_{r}(u)}{n^{1-\epsilon}} j_{\tau_{0}}(u)
\]
and
\[
   Z_{12}(s,\alpha) = \prod_{p^{\lambda} \mid \mid v,p \le c_{3}} (c_{5}d_{r}(p^{\lambda}))
   \prod_{p^{\lambda} \mid \mid v, p > c_{3}} d_{r}(p^{\lambda})(1 + O(p^{-1+2 \epsilon}))
   \ll d_{r}(v) j_{\tau_{0}}(v) \ .
\]
We conclude that if $\mathrm{Re}(s) \ge -\epsilon$ and $|\alpha| \le \epsilon$ then
$
   |Z_{1}(s,\alpha)| 
   \ll d_{r}(m) j_{\tau_{0}}(m) n^{\epsilon-1} \ .
$
This completes our calculation of~(\ref{eq:z1s}).  The lemma will be thus
completed once the residue is computed.
\subsection{The residue computation}
We decompose
\begin{equation}
  Z(s,\alpha)x^{s}s^{-1} = \zeta(1+s-i\alpha)^{-r} 
  Z_{1}(s,\alpha)Z_{3}(s,\alpha) x^{s}
  \zeta(1+s)^{2r}s^{-1} \ .
  \label{eq:prod}
\end{equation}
We now compute the Laurent expansion of each factor.  We have 
\begin{equation}
\begin{split}
  \zeta^{2r}(1+s)s^{-1} & = s^{-2r-1} (1 + a_{1}s + a_{2} s^{2} + \cdots) \ , \\
  x^{s} & = 1 + (\log x) s + (\log x)^2 s^{2}/2! + \cdots   \ , \\
  \zeta(1+s-i \alpha)^{-r} & = f(-i \alpha) 
  + f^{'}(-i \alpha)s + f^{(2)}(-i \alpha) s^{2}/2! + \cdots\\
   \label{eq:laurent}
\end{split}
\end{equation}
where we put $f(z)= \zeta(1+z)^{-r}$.  Note that 
a simple calculation yields
\begin{equation}
  f^{(j)}(-i \alpha) = \left\{ \begin{array}{cl}
                  r(r-1) \cdots (r-(j-1))(-i \alpha)^{r-j} 
                  + O(|\alpha|^{r-j+1})  &  \mbox{$0 \le j \le r$} \\
                  c_{j} + O(|\alpha|) &  \mbox{$j \ge r+1$}  \\
                  \end{array}
                  \right.
\end{equation}
and $c_{j} \in \mathbb{R}$.  Next note that $Z_3(s,\alpha)$ has 
an absolutely convergent power series in  $\mathrm{Re}(s) > -1/2$, $|\alpha| \le cL^{-1}$. It follows that $Z_{3}(0,\alpha) = Z_{3}(0,0) + O(|\alpha|) = 1 + O(|\alpha|)$
and $Z^{(j)}_{3}(0,\alpha) \ll 1$
for $j \ge 0$.  Combining these facts yields
\begin{equation}
  Z_{3}(s,\alpha) = (1 + O(|\alpha|)) + O(1)s + O(1)s^{2} + \cdots \ .
  \label{eq:z3taylor}
\end{equation}
We now compute the Taylor expansion of $Z_{1}(s,\alpha)$.  
Since $k_{p}(0)=0$ it follows from
~(\ref{eq:alps}),~(\ref{eq:beps}),~(\ref{eq:hps}) , and~(\ref{eq:z11z12}) that
\begin{equation}
  Z_{1}(s,0)   = \frac{\sigma_{r}(m,s+1)}{n} \ .
\end{equation}
By Cauchy's integral formula with a circle of radius $\epsilon/2$, we 
establish a bound for $Z_{1}^{(j)}(0,\alpha)$:
\begin{equation}
     Z_{1}^{(j)}(0,\alpha)
     = \frac{1}{2 \pi i} \int_{|w-\alpha|=\epsilon/2}
     \frac{Z_{1}(0,w) dw}{(w-\alpha)^{j+1}}
     \ll \left( \frac{2}{\epsilon} \right)^{j+1}
     \frac{d_{r}(m)j_{\tau_0}(m)}{n^{1-\epsilon}}
     \label{eq:z1j0s}
\end{equation}
by~(\ref{eq:z1s}).  
By the Taylor series expansion and~(\ref{eq:z1j0s}) it follows that
\begin{equation}
  Z_{1}(0,\alpha) =  
  \frac{\sigma_{r}(m)}{n} + O \left( \frac{d_{r}(m)j_{\tau_{0}}(m)}{n^{1-\epsilon}}|\alpha|  \right)
     \label{eq:z10al}
\end{equation}
since $Z_1(0,0) = \sigma_{r}(m)/n$.
Combining~(\ref{eq:z1j0s}) and~(\ref{eq:z10al}) we obtain 
\begin{equation}
  Z_{1}(s,\alpha) =
  \left( \frac{\sigma_{r}(m)}{n} + O \left( \frac{d_{r}(m)j_{\tau_{0}}(m)}{n^{1-\epsilon}}|\alpha|  \right)   \right)
  +  \sum_{j=0}^{\infty} O(d_{r}(m)j_{\tau_{0}}(m)n^{\epsilon-1}) s^{j}/j! \ .
  \label{eq:z1taylor}
\end{equation}
We are now in a position to compute the residue.
It follows from~(\ref{eq:prod}),~(\ref{eq:laurent}),~(\ref{eq:z3taylor}),
and~(\ref{eq:z1taylor}) that the residue at $s=0$ is
\begin{equation}
  res = \sum_{u_{1}+u_{2}+u_{3}+u_{4}+u_5=2r} \frac{l^{u_{1}}f^{(u_{2})}(-i \alpha)Z_{1}^{(u_{3})}(0,\alpha)
  Z_{3}^{(u_{4})}(0,\alpha) a_{u_5}}{u_1! u_{2}! u_{3}! u_{4}!}  \ . 
  \label{eq:res0}
\end{equation}
We first show that those terms with $u_5 \ge 1$ contribute a smaller 
amount.  Since $|f^{(u_2)}(-i \alpha)| \ll |\alpha|^{r-u_2}$ for $0 \le u_2 \le r$
and $|f^{(u_2)}(-i \alpha)|  \ll_{r} 1$ for $r+1 \le u_2 \le 2r$ it follows that
the terms with $u_5 \ge 1$ contribute
\begin{equation}
\begin{split}
    & \ll_{r} \frac{d_{r}(m)j_{\tau}(m)}{n^{1-\epsilon}}
   \sum_{u_1+u_2 \le 2r-1} \frac{l^{u_1}| f^{(u_2)}(-i \alpha)|}{u_1! u_2!}  \\
   & 
   \ll  \frac{d_{r}(m)j_{\tau}(m)}{n^{1-\epsilon}}
   \left(
   \sum_{{\begin{substack}{u_{1}+u_{2} \le 2r-1
         \\ 0 \le  u_{2} \le r}\end{substack}}} l^{u_1} |\alpha|^{r-u_2}
    +  \sum_{{\begin{substack}{u_{1}+u_{2} \le 2r-1
         \\ 0 \le  u_{2} \ge r+1}\end{substack}}} l^{u_1} 
   \right)  \\
   & \ll 
   \frac{d_{r}(m)j_{\tau}(m)}{n^{1-\epsilon}}
  (l^{r-1} + l^{r-2})  \ . 
   \nonumber 
\end{split}
\end{equation}
We deduce that 
\begin{equation}
  res = \sum_{u_{1}+u_{2}+u_{3}+u_{4}+u_5=2r} \frac{l^{u_{1}}f^{(u_{2})}(-i \alpha)Z_{1}^{(u_{3})}(0,\alpha)
  Z_{3}^{(u_{4})}(0,\alpha)}{u_1! u_{2}! u_{3}! u_{4}!} 
  + 
   O \left( \frac{d_{r}(m) j_{\tau_{0}}(m) l^{r-1}}{n^{1-\epsilon}}  \right)
   \ . 
  \label{eq:res}
\end{equation}
The contribution from those terms in~(\ref{eq:res}) satisfying 
$u_1+u_2=2r, u_{2} \le r$ is
\begin{equation}
\begin{split}
  & \left( \sum_{u_{1}+u_{2}=2r, u_2 \le r} \frac{l^{u_1}}{u_1!} \frac{f^{(u_2)}(-i \alpha)}{u_2!} \right)
  \left( \frac{\sigma_{r}(m)}{n} + O \left( \frac{d_{r}(m)j_{\tau_{0}}(v)}{n^{1-\epsilon}}|\alpha|  \right) \right)
  ( 1 + O(|\alpha|)) \ . \\
  & = \left( \sum_{u_{2} \le r} \frac{l^{2r-u_{2}}}{(2r-u_{2})!} ( \binom{r}{u_{2}} (-i \alpha)^{r-u_{2}}
  + O(|\alpha|^{r-u_{2}+1}) )\right)
   \left( \frac{\sigma_{r}(m)}{n} + O \left( \frac{d_{r}(m)j_{\tau_{0}}(v)}{n^{1-\epsilon}}|\alpha|  \right) \right) \\
   & = \frac{\sigma_{r}(m)}{n} l^{r}  \sum_{a=0}^{r} \binom{r}{a} \frac{1}{(r+a)!} (-i \alpha l)^{a} +
   O \left( \frac{d_{r}(m) j_{\tau_{0}}(m)}{n^{1-\epsilon}}  \right) \ .
  \nonumber
\end{split}
\end{equation}
Those terms in~(\ref{eq:res}) with $u_{1} \le r-1$ contribute
\begin{equation}
  \frac{d_{r}(m) j_{\tau_{0}}(m)}{n^{1-\epsilon}}
   \sum_{{\begin{substack}{u_{1}+u_{2}+u_{3}+u_{4}=2r
         \\ u_{1} \le r-1}\end{substack}}} l^{u_1}|f^{(u_{2})}(- i \alpha)|  \ll
   \frac{d_{r}(m) j_{\tau_{0}}(m)}{n^{1-\epsilon}} l^{r-1}
   \nonumber
\end{equation}
since $|\alpha| \le cL^{-1} \ll 1$ and the remaining terms in~(\ref{eq:res}) 
are
\begin{equation}
   \ll \frac{d_{r}(m) j_{\tau_{0}}(m)}{n^{1-\epsilon}}
   \sum_{{\begin{substack}{u_1+u_2+u_3+u_4=2r
         \\ u_1+u_2 \le 2r-1 , r \le  u_1 \le 2r-1}\end{substack}}} l^{u_1} |\alpha|^{u_{1}+1-r}
   \ll     \frac{d_{r}(m) j_{\tau_{0}}(m)}{n^{1-\epsilon}} l^{r-1} \ .
   \nonumber
\end{equation}
We thus conclude that 
\begin{equation}
   res = \frac{\sigma_{r}(m)}{n} l^{r}  \sum_{a=0}^{r} \binom{r}{a} \frac{1}{(r+a)!} (-i \alpha l)^{a} +
   O \left( \frac{d_{r}(m) j_{\tau_{0}}(m) l^{r-1}}{n^{1-\epsilon}}  \right)
   \label{eq:resf}
\end{equation}
and the lemma follows from 
~(\ref{eq:main}) and~(\ref{eq:resf}). 

We have the Taylor series expansion
$$
   \mathcal{R}_{k}(1+i \alpha)
  = \rk_{k}(1) + \rk^{'}_{k}(1)(i \alpha) +
  \rk_{k}^{(2)}(1) (i \alpha)^{2}/2 + \cdots \ .
$$
We denote the truncated Taylor series expansion $\mathcal{T}_{k;N}(\alpha) = \sum_{j=0}^{N} \rk_{k}^{(j)}(1)(i
\alpha)^{j}/j!$.
\begin{lemma} We
have for $l = \log x$, $|\alpha| \ll (\log x)^{-1}$, and $\tau_{0} = 1/3$
\begin{equation}
   \sum_{k \le x} d(mk) \frac{\mathcal{T}_{nk;r}(\alpha)}{\phi(nk)}
   = \frac{\sigma_{r}(m)}{n} l^{r}
   \sum_{j=0}^{r} \binom{r}{j} \frac{1}{(r+j)!} (-i \alpha \ l)^{j} +
   O \left( \frac{d_{r}(m)j_{\tau_{0}}(v)l^{r-1}}{n^{1-\epsilon}}  \right)
   .
\end{equation}
\end{lemma}
\noindent {\it Proof}.  We begin by noting that it suffices to prove
\begin{equation}
  \sum_{k \le x} \frac{d_{r}(mk) \rk_{nk}^{(j)}(1)}{\phi(nk)}
  = (-1)^{j} \frac{\sigma_{r}(m)}{n} l^{r+j} \binom{r}{j} \frac{j!}{(r+j)!}
  + O \left(\frac{d_{r}(m) j_{\tau_{0}}(m)l^{r+j-1}}{n}  \right) \ .
  \label{eq:drrn}
\end{equation}
This is since if we multiply~(\ref{eq:drrn})
by $(i \alpha)^{j}/j!$ and sum $j=$ to $r$ we obtain the result. The
Dirichlet series generating function for the sum in question is
\begin{equation}
  i^{-j} \sum_{k=1}^{\infty} \frac{d_{r}(mk)}{\phi(nk) k^{s}}
   \left. \frac{d^{j}}{d \alpha^{j}} \rk_{nk}(1+i \alpha) \right|_{\alpha = 0}
  = i^{-j} \left. \frac{d^{j}}{d \alpha^{j}} Z(s,\alpha) \right|_{\alpha=0} \ .
\end{equation}
By Perron's formula it follows that the sum in question is
\begin{equation}
  \frac{i^{-j}}{2 \pi i} \int_{c-iU}^{c+iU}  \left. \frac{d^{j}}{d \alpha^{j}} Z(s,\alpha) \right|_{\alpha=0}
  \frac{x^{s}}{s} \, ds + O \left( \frac{d_{r}(m)}{n^{1-\epsilon}}
  \left(\frac{(\log x)^{r}}{U} + 1   \right)     \right) \ .
\end{equation}
where $c=(\log x)^{-1}$.  As in Lemma 13 equations~(\ref{eq:perr}),
~(\ref{eq:deform}),  we want to deform the
contour $[c-iU,c+iU]$ to $\Gamma(U)$ and then pick up the residue at
$s=0$.  As this calculation is analogous to the preceding lemma we omit
the details.  This procedure yields
\begin{equation}
\sum_{k \le x} \frac{d_{r}(mk) \rk_{nk}^{(j)}(1)}{\phi(nk)}
  = i^{-j} 
   {\begin{substack} { res \\ s=0} \end{substack}} 
   \left(
    \left. \frac{d^{j}}{d \alpha^{j}} Z(s,\alpha) \right|_{\alpha=0}
  \frac{x^{s}}{s} \right)
   +
   O \left( \frac{d_{r}(m)j_{\tau_{0}}(v)l^{r-1}}{n^{1-\epsilon}}  \right)
   \ . 
   \label{eq:res14}
\end{equation}
Recall that $Z(s,\alpha) = Z_{1}(s,\alpha)Z_{2}(s,\alpha)Z_{3}(s,\alpha)$ where
\begin{equation}
  Z_{1}(s,0) = \frac{\sigma_{r}(m,s+1)}{n} \ , \
  Z_{2}(s,\alpha) = \frac{\zeta^{2r}(1+s)}{\zeta^{r}(1+s-i \alpha)} \ ,
  \ Z_{3}^{(j)}(0,0) \ll 1
\end{equation}
for all $j \ge 0$.  By the product rule we have
\begin{equation}
   \left.  \frac{d^{j}}{d \alpha^{j}} Z(s,\alpha)  \right|_{\alpha=0}
   = \sum_{u_{1}+u_{2}+u_{3}=j} \binom{j}{u_{1},u_{2},u_{3}}
   Z_{1}^{(u_{1})}(s,0) Z_{2}^{(u_{2})}(s,0) Z_{3}^{(u_{3})}(s,0) \ .
   \label{eq:res14a}
\end{equation}
Thus we need to compute 
\begin{equation}
   \mathrm{res}_{s=0} \left(
Z_{1}^{(u_{1})}(s,0) Z_{2}^{(u_{2})}(s,0) Z_{3}^{(u_{3})}(s,0)
 x^{s} s^{-1} \right)
   \label{eq:res14b}
 \end{equation}
 for all $u_1+u_2+u_3=j$. 
In fact, it turns out that the main 
term arises from those triples $(u_{1},u_{2},u_{3})=(0,j,0)$. 
We now compute the residue arising from these terms. We
have the Laurent expansions,
\begin{equation}
\begin{split}
   Z_{1}(s,0)  & = \frac{\sigma_{r}(m)}{n} + \frac{\sigma_{r}^{(1)}(m,1)}{n} s + \cdots  \ ,  \\
   Z_{2}^{(j)}(s,0)  & =  \frac{r(r-1) \cdots (r-(j-1)) (-i)^{j}}{s^{r+j}} +  \frac{c_{1}}{s^{r+j-1}}  + \cdots   \ ,  \\
   Z_{3}(s,0) & = 1 + d_{1} s + \cdots \ .  
   \nonumber
\end{split}
\end{equation}
We further remark that by Cauchy's integral formula we may establish
$\sigma_{r}^{(k)}(m,1) \ll d_{r}(m) j_{\tau_{0}}(m)$
for some $\tau_{0} >0$.  These terms contribute
\begin{equation}
\begin{split}
    &  {\begin{substack} { res \\ s=0} \end{substack}} Z_{1}(s,0)Z_{2}^{(j)}(s,0)Z_{3}(s,0)x^{s} s^{-1} = \\
   & \frac{\sigma_{r}(m) r(r-1) \cdots (r-(j-1)) l^{r+j} (-i)^{j}}{n (r+j)!} +
   O \left(
   \frac{d_{r}(m) j_{\tau_{0}}(m) l^{r+j-1}}{n}
   \right) \ .
   \label{eq:res14c}
\end{split}
\end{equation}
A similar calculation shows that for those triples $(u_{1},u_{2},u_{3})$ such that $u_{2} \le j-1$ then
\begin{equation}
{\begin{substack} { res \\ s=0} \end{substack}} Z_{1}^{(u_{1})}(s,0)Z_{2}^{(u_{2})}(s,0)
                  Z_{3}^{(u_{3})}(s,0)x^{s}s^{-1}
  \ll \frac{d_{r}(m) j_{\tau_{0}}(m)(\log x)^{r+j-1}}{n} \ .
  \label{eq:res14d}
\end{equation}
The lemma now follows by combining
~(\ref{eq:res14}),~(\ref{eq:res14a}),~(\ref{eq:res14b}),~(\ref{eq:res14c}), and~(\ref{eq:res14d}). 

We deduce the following corollary to Lemmas 13 and 14:
\begin{lemma} 
    \begin{equation}
   \sum_{k \le x} d_{r}(mk) \left(f(nk)-  \frac{\mathcal{T}_{nk,r}(\alpha)}{\phi(nk)} \right)
   \ll |\alpha|^{r+1} L^{2r} \frac{d_{r}(m)j_{\tau_{0}}(m)}{n^{1-\epsilon}}  \ .
   \label{eq:coro}
\end{equation}
\end{lemma}
\noindent {\it Proof}.  Note that
\[
   f(nk) = \frac{\mathcal{T}_{nk,r}(\alpha)}{\phi(nk)} + \alpha^{r+1}
   g(\alpha;nk)
\]
where $g$ is entire in $\alpha$.  Moreover, it follows that
\begin{equation}
   \sum_{k \le x} d_{r}(mk) \left(f(nk)-\frac{\mathcal{T}_{nk,r}(\alpha)}{\phi(nk)} \right)
   = \alpha^{r+1} g^{*}(\alpha;n,x)
   \label{eq:defg}
\end{equation} where $g^{*}$ entire in $\alpha$.  Combining 
Lemmas 13 and 14 we deduce that 
\[
  \max_{|\alpha| \le cL^{-1}} |\alpha^{r+1} g^{*}(\alpha;n,x)|
  \ll \frac{d_{r}(m)j_{\tau_{0}}(m) L^{r-1}}{n^{1-\epsilon}}
\]
and hence by the maximum modulus principle
\begin{equation}
    \max_{|\alpha| \le cL^{-1}} |g^{*}(\alpha;n,x)|
  \ll \frac{d_{r}(m)j_{\tau_{0}}(m)L^{2r}}{n^{1-\epsilon}} \ .
  \label{eq:mmp}
\end{equation}
Hence,~(\ref{eq:defg}) and~(\ref{eq:mmp}) imply the statement of the 
lemma.

\section{Proof of Theorem 2}

\subsection{Initial manipulations}
In this section we apply the lemmas to manipulate $I$ into a suitable form for evaluation.  Recall that we had~(\ref{eq:Iexp})
\begin{equation}
  I = \sum_{k \le y} \frac{d_{r}(k)P([k]_{y})}{k} \sum_{j \le \frac{kT}{2
  \pi}} b(j) e (-j/k) + O_{\epsilon}(y
  T^{\frac{1}{2}+ \epsilon})  \ . 
  \label{eq:e51}
\end{equation}
By Perron's formula with $c=1+L^{-1}$ the inner sum is
\[
  \sum_{j \le \frac{kT}{2
  \pi}} b(j) e (-j/k)
  = \frac{1}{2 \pi i} \int_{c-iT}^{c+iT} \Q^{*}(s,\alpha,k) \left(
  \frac{kT}{2 \pi} \right)^{s} \frac{ds}{s} + O(kT^{\epsilon}) 
\]
where $\Q^{*}(s,\alpha,k) = \sum_{j=1}^{\infty} b(j)j^{-s}e (-j/k)$. 
Pulling the contour left to $c_{0} = 1/2 + L^{-1}$ we obtain
\begin{equation}
\begin{split}
   \sum_{j \le \frac{kT}{2
  \pi}} b(j) e (-j/k)  & =
  R_{1} + R_{1+i \alpha}  \\
  &  + \frac{1}{ 2 \pi i} \left(
   \int_{c-iT}^{c_{0}-iT} + \int_{c_{0}-iT}^{c_{0}+iT} +
   \int_{c_{0}+iT}^{c+iT}
   \right) \Q^{*}(s,\alpha,k)
   \left( \frac{kT}{2 \pi }\right)^{s} \frac{ds}{s}
   \label{eq:e52}
\end{split}
\end{equation}
where $R_{u}$ is the residue at $s=u$. By Lemma 8 the left and horizontal edges 
contribute $yT^{1/2+\epsilon}$.
Moreover by~(\ref{eq:bj})  it follows that
\begin{equation}
   \Q^{*}(s,\alpha,k) = \sum_{h \le y} \frac{d_{r}(h) P([h]_{y}) \Q(s,\alpha,h/k)}{h^{s}}
\end{equation}
where $\Q(s,\alpha,h/k)$ is defined by~(\ref{eq:qsa}).  
We will now invoke Lemma 5, however we require that $h,k$ be relatively prime. 
Therefore we set $\frac{h}{k}= \frac{H}{K}$ where
$H  = h/(h,k)$, $K = k/(h,k)$, 
and $(H,K)=1$.  We deduce
\[
     R_{1} = \sum_{h \le y} d_{r}(h)P([h]_{y}) \  {\begin{substack} { res \\ s=1} \end{substack}}  \left(\Q(s,\alpha,H/K)
    \left( \frac{TK}{2 \pi H} \right)^{s} s^{-1}
       \right) \ .
\]
By an application of Lemma 5$(i)$ this is
\begin{equation}
\begin{split}
    & R_{1} = K \sum_{h \le y} d_{r}(h) P([h]_{y}) \   {\begin{substack} { res \\ s=1} \end{substack}}  \left(
    \zeta^{2}(s) \left( \frac{\zeta^{'}}{\zeta}(s-i \alpha)
    - \mathcal{G}(s,\alpha,K)  \right)
    \left( \frac{T}{2 \pi HK} \right)^{s} s^{-1}
    \right) \\
    & = \frac{T}{2 \pi}
    \sum_{h \le y} \frac{d_{r}(h) P([h]_{y})}{H}  \\
    & \cdot  \left(
    \left( (\zeta^{'}/\zeta)(\overline{\tau}) - \mathcal{G}(1, \alpha, K)
    \right)
    \log \left( \frac{Te^{2 \gamma -1}}{2 \pi H K} \right)
    + \left(
    (\zeta^{'}/\zeta)^{'}(\overline{\tau}) -
    \mathcal{G}^{'}(1, \alpha, K)
    \right)
    \right)
    \label{eq:e53}
\end{split}
\end{equation}
where we put $\tau = 1+i \alpha$.  Likewise Lemma 5$(ii)$ implies
\begin{equation}
\begin{split}
    R_{1+i \alpha} & =
    \sum_{h \le y} d_{r}(h) P([h]_{y})  \  {\begin{substack} { res \\ s=\tau} \end{substack}}  \left(\Q(s,\alpha,H/K)
    \left( \frac{TK}{2 \pi H} \right)^{s} s^{-1} \right) \\
     & = - \frac{T}{2 \pi} \frac{\zeta^{2}(\tau)}{\tau}
       \sum_{h \le y} \frac{d_{r}(h)P([h]_{y}) }{H} \left( \frac{T}{2 \pi H} \right)^{i \alpha}
       \frac{K\rk_{K}(\tau)}{\phi(K)} \ .
   \label{eq:e54}
\end{split}
\end{equation}
Combining~(\ref{eq:e51}), ~(\ref{eq:e52}),~(\ref{eq:e53}), and
~(\ref{eq:e54}) we deduce
\[
   I = \frac{T}{2 \pi}
           \sum_{h,k \le y} \frac{d_{r}(h)d_{r}(k)P([h]_{y}) P([k]_{y}) (h,k)}{hk}
           \left( \log \frac{Te^{2 \gamma-1}}{2 \pi HK}
           \left(
           (\zeta^{'}/\zeta)(\tau)-\mathcal{G}(1,\alpha,K)\right)
           \right.
\]
\[
     \left.
     + \left( \zeta^{'}/\zeta\right)^{'}(\tau)
     - \mathcal{G}^{'}(1,\alpha,K)
     - \frac{\zeta^{2}(\tau)}{\tau}
       \left( \frac{T}{2 \pi H} \right)^{i\alpha}
       \frac{K\rk_{K}(\tau)}{\phi(K)} \right) + O(yT^{1/2+\epsilon})
\]
where $\mathcal{G}(s,\alpha,K)$ is defined by~(\ref{eq:gsak}).  We may write for $j=0,1$
$\mathcal{G}^{(j)}(1,\alpha,K) = \sum_{p \mid K} p^{i\alpha} \log^{j+1} p + O(C_{j}(K))$.  By Lemma 9, the
$O(C_{j}(K))$ terms contribute $O(TL^{(r+1)^2})$. Whence
\[
   I = \frac{T}{2 \pi}
           \sum_{h,k \le y} \frac{d_{r}(h)d_{r}(k)P([h]_{y}) P([k]_{y}) (h,k)}{hk}
           \left( \log \frac{Te^{2 \gamma-1}}{2 \pi HK}
           \left(
           (\zeta^{'}/\zeta)(\overline{\tau})
           - \sum_{p \mid K} p^{i\alpha} \log p \right)
           \right.
\]
\[
     \left.
     + (\zeta^{'}/\zeta)^{'}(\overline{\tau})
     - \sum_{p \mid K} p^{i\alpha}
     \log^{2} p
     - \frac{\zeta^{2}(\tau)}{\tau}
       \left( \frac{T}{2 \pi H} \right)^{i\alpha}
       \frac{K \rk_{K}(\tau)}{\phi(K)} \right)  + O(y T^{1/2+\epsilon}) 
\]
where $z=1+i \alpha$.  Insertion of the identity
\[
  f((h,k)) = \sum_{{\begin{substack}{m \mid h \\
                             m \mid k}
            \end{substack}}} \sum_{n \mid m} \mu(n) f
  \left( \frac{n}{m} \right)
\]
produces
\begin{equation}
\begin{split}
  I =  \frac{T}{2 \pi} &
           \sum_{h,k \le y} \frac{d_{r}(h)P([h]_{y})d_{r}(k)P([k]_{k})}{hk} \sum_{
           {\begin{substack}{m \mid h \\
                             m \mid k}
            \end{substack}}}
            m \sum_{n \mid m} \frac{\mu(n)}{n}   \\
    & \cdot \left( \log \frac{Te^{2 \gamma-1}m^{2}}{2 \pi hkn^{2}}
           \left(
           (\zeta^{'}/\zeta)(\overline{\tau})- \sum_{p \mid \frac{nk}{m}} p^{i\alpha}
          \log p \right)
           \right. + \left( \zeta^{'}/\zeta \right)^{'}(\overline{\tau})
    \\
    &
     \left.
     - \sum_{p \mid \frac{nk}{m}} p^{i\alpha}
     \log^{2} p
     - \frac{\zeta^{2}(\tau)}{\tau}
       \left( \frac{Tm}{2 \pi nh} \right)^{i\alpha} \left(
       \frac{nk}{m} \right)
       \frac{\rk_{\frac{nk}{m}}(\tau)}{\phi(\frac{nk}{m})} \right)
       + O(y T^{1/2+\epsilon}) \ . 
       \nonumber
\end{split}
\end{equation}
Changing summation order and making the variable changes $h \to hm$ and $k \to km$ yields
\begin{equation}
\begin{split}
  I = \frac{T}{2 \pi} & \sum_{m \le y} \frac{1}{m} \sum_{n \mid m}
  \frac{\mu(n)}{n}
  \sum_{h,k \le \frac{y}{m}} \frac{d_{r}(mh)P([mh]_{y})d_{r}(mk) P([mk]_{y})}{hk} \\
  &     \cdot   \left( \log \frac{Te^{2 \gamma-1}}{2 \pi hkn^{2}}
           \left(
           (\zeta^{'}/\zeta)(\overline{\tau})- \sum_{p \mid nk} p^{i\alpha}
          \log p \right)
           \right. + \left( \zeta^{'}/\zeta \right)^{'}(\overline{\tau})\\
  &
     \left.
     - \sum_{p \mid nk} p^{i\alpha}
     \log^{2} p
     - \frac{\zeta^{2}(\tau)}{\tau}
       \left( \frac{T}{2 \pi nh} \right)^{i\alpha}
       \frac{nk \rk_{nk}(\tau)}{\phi(nk)} \right) 
       + O(y T^{1/2+\epsilon}) \ .
       \nonumber
\end{split}
\end{equation}
Rearrange this as $I = I_{1} + I_{2}+ O(yT^{1/2+\epsilon})$ where
\begin{equation}
\begin{split}
   I_{1} =  \frac{T}{2 \pi} &
   \sum_{m \le y} \frac{1}{m}
   \sum_{n \mid m} \frac{\mu(n)}{n}
   \sum_{h,k \le \frac{y}{m}} \frac{d_{r}(mh)P([mh]_{y})d_{r}(mk)P([mk]_{y})}{hk}  \\
   & \cdot \left(
   -\log \frac{Te^{2 \gamma-1}}{2 \pi hkn^{2}} \sum_{p \mid nk}
   p^{i\alpha} \log p - \sum_{p \mid nk} p^{i \alpha} \log^{2} p
   \right)
   \nonumber
\end{split}
\end{equation}
and
\begin{equation}
\begin{split}
   I_{2} = \frac{T}{2 \pi} &
   \sum_{m \le y} \frac{1}{m}
   \sum_{n \mid m} \frac{\mu(n)}{n}
   \sum_{h,k \le \frac{y}{m}} \frac{d_{r}(mh)P([mk]_{y})d_{r}(mk)P([mk]_{y})}{hk}  \\
   & \cdot \left( \log \left( \frac{Te^{2 \gamma-1}}{2 \pi hkn^{2}} \right)
   \frac{\zeta^{'}}{\zeta} (\overline{\tau})
   + \left(  \frac{\zeta^{'}}{\zeta} \right)^{'}(\overline{\tau})
   \right.
    \left.
   - \frac{\zeta^{2}(\tau)}{\tau}
   \left( \frac{T}{2 \pi nh} \right)^{i \alpha}
   \frac{nk \rk_{nk}(\tau)}{\phi(nk)} \right) \ .
   \label{eq:I2}
\end{split}
\end{equation}
The first sum is
\begin{equation}
\begin{split}
   I_{1} & = \frac{T}{2 \pi}
   \sum_{m \le y} \frac{1}{m}
   \sum_{n \mid m} \frac{\mu(n)}{n}
   \sum_{h,k \le \frac{y}{m}} \frac{d_{r}(mh)P([mh]_{y})d_{r}(mk)P([mk]_{y})}{hk} \\
   &  \cdot \left(
   -\log \frac{T}{2 \pi hk} \sum_{p \mid k}
   p^{i\alpha} \log p - \sum_{p \mid k} p^{i \alpha} \log^{2} p
   + O(L \log n) \right) \ .
    \nonumber
\end{split}
\end{equation}
A calculation shows that the $O(L \log n)$ contributes $O(TL^{(r+1)^2})$. Since $\phi(m)m^{-1} = \sum_{n \mid m}
\mu(n) n^{-1}$ we deduce that
\begin{equation}
\begin{split}
  I_{1} = \frac{T}{2 \pi}
   & \sum_{m \le y} \frac{\phi(m)}{m^{2}}
   \sum_{h,k \le \frac{y}{m}} \frac{d_{r}(mh)P([mh]_{y})d_{r}(mk)P([mk]_{y})}{hk} \\
   & \cdot \left(
   -\log \frac{T}{2 \pi hk} \sum_{p \mid k}
   p^{i\alpha} \log p - \sum_{p \mid k} p^{i \alpha} \log^{2} p
   \right) + O(TL^{(r+1)^2}) \ .
   \label{eq:I1b}
\end{split}
\end{equation}
This puts $I_{1}$ in a suitable form to be evaluated by the lemmas. We now simplify $I_{2}$ by substituting the
Laurent expansions
\begin{equation}
\begin{split}
    (\zeta^{'}/\zeta)(\overline{\tau})
   & = (i \alpha)^{-1} + O(1) \ , \nonumber  \\
   \left(\zeta^{'}/\zeta \right)^{'} (\overline{\tau}) & = (i \alpha)^{-2} +
   O(1) \ , \nonumber  \\
   \zeta^{2}(\tau) \tau^{-1}
   & = (i \alpha)^{-2} +  (2 \gamma-1)(i \alpha)^{-1} + O(1) \
   \nonumber
\end{split}
\end{equation}
in~(\ref{eq:I2}).  The $O(1)$ terms of these Laurent expansions contribute
\[
   TL \sum_{m \le y} \frac{d_{r}(m)^{2}}{m} \sum_{n \mid m} \frac{1}{n}
  \sum_{h,k \le \frac{y}{m}} \frac{d_{r}(h)d_{r}(k)}{hk}
   \ll TL^{(r+1)^{2}}
\]
by~(\ref{eq:ram}) and
\begin{equation}
\begin{split}
  & T \sum_{m \le y} \frac{d_{r}(m)}{m} \sum_{n \mid m} 1
  \left( \sum_{h \le \frac{y}{m}} \frac{d_{r}(h)}{h} \right)
  \left| \sum_{k \le \frac{y}{m}}  d_{r}(mk)f(nk) \right| \\
  & \ll TL^{r} \sum_{m \le y} \frac{d_{r}(m)}{m}
      \sum_{n \mid m} \left( \frac{d_{r}(m)\sigma_{r}(m)L^{r}}{n^{1-\epsilon}}
      \right) \ll TL^{(r+1)^{2}-1}
  \nonumber
\end{split}
\end{equation}
by~(\ref{eq:ram}) and Lemma 13. Thus we deduce
\begin{equation}
\begin{split}
   I_{2} = \frac{T}{2 \pi}
   \sum_{m \le y} \frac{1}{m} & \sum_{n \mid m} \frac{\mu(n)}{n}
   \sum_{h,k \le \frac{y}{m}} \frac{d_{r}(mh)P([mh]_{y})d_{r}(mk)P([mk]_{y})}{hk} \\
   & \cdot \left(
   \frac{1 + i \alpha \log \frac{T}{2 \pi hk n^{2}} -
   \left( \frac{T}{2 \pi h n} \right)^{i \alpha}
   \frac{nk \rk_{nk}(\tau)}{\phi(nk)}}{(i \alpha)^{2}}
   \right)
   \nonumber
\end{split}
\end{equation}
plus an error term $O(TL^{(r+1)^{2}})$. In the above formula we replace $\frac{\rk_{nk}(\tau)}{\phi(nk)}$ by
$\frac{\mathcal{T}_{nk;r}(\alpha)}{\phi(nk)}$
and by~(\ref{eq:coro}) this introduces an error of
\begin{align}
  & |\alpha|^{-2} T L^{r}  \sum_{m \le y} \frac{d_{r}(m)}{m} \sum_{n \mid m}
   \left|  \sum_{k \le \frac{y}{m}} d_{r}(mk)
    \left(
   \frac{\mathcal{R}_{nk}(\tau)}{\phi(nk)} - \frac{\mathcal{T}_{nk;r}(\alpha)}{\phi(nk)}
   \right) \right|  \nonumber \\
   &  \ll |\alpha|^{-2} TL^{r}  \sum_{m \le y} \frac{d_{r}(m)}{m} \sum_{n \mid m}
   \left( |\alpha|^{r+1} L^{2r}
    \frac{d_{r}(m)j_{\tau_0}(m)}{n^{1-\epsilon}} \right) \ll TL^{(r+1)^2}  \ .
    \nonumber
\end{align}
Therefore we have
\begin{equation}
\begin{split}
   I_{2} & = \frac{T}{2 \pi}
   \sum_{m \le y} \frac{1}{m} \sum_{n \mid m} \frac{\mu(n)}{n}
   \sum_{h,k \le \frac{y}{m}} \frac{d_{r}(mh)P([mh]_{y})d_{r}(mk)P([mk]_{y})}{hk}  \\
   & \cdot \left(
   \frac{1 + i \alpha \log \frac{T}{2 \pi hk n^{2}} -
   \left( \frac{T}{2 \pi h n} \right)^{i \alpha}
   \frac{nk \mathcal{T}_{nk;r}(\alpha)}{\phi(nk)}}{(i \alpha)^2}
   \right)  + O(TL^{(r+1)^2})  \ .
   \nonumber
\end{split}
\end{equation}
A calculation shows that $\rk_{k}(1) = \phi(k)/k$, $\rk_{k}^{'}(1) = - \phi(k) \log k/k$ and thus it follows that
\begin{equation}
  \frac{\mathcal{T}_{nk;r}(\alpha)}{\phi(nk)}  =  \frac{1}{nk} \left( 1 - \log(nk) (i \alpha) \right)
  + \sum_{j=2}^{r} \frac{\rk_{nk}^{(j)}(1)(i \alpha)^{j}}{\phi(nk) j!} \ .
  \nonumber
\end{equation}
We further decompose $I_{2} = I_{21} + I_{22} + O(TL^{(r+1)^2})$ where
\begin{equation}
\begin{split}
  I_{21} = \frac{T}{2 \pi}
   & \sum_{m \le y} \frac{1}{m} \sum_{n \mid m} \frac{\mu(n)}{n} \sum_{h,k \le \frac{y}{m}}
   \frac{d_{r}(mh)P([mh]_{y})d_{r}(mk)P([mk]_{y})}{hk}  \\
   & \cdot \left(
   \frac{1 + i \alpha \log \frac{T}{2 \pi hk n^{2}} -
   \left( \frac{T}{2 \pi h n} \right)^{i \alpha}
   (1-(i \alpha) \log (nk))}{(i \alpha)^2}
   \right)
\end{split}
\end{equation}
\begin{equation}
\begin{split}
   \mathrm{and} \
   I_{22}  & = - \frac{T}{2 \pi} \sum_{j=2}^{r} \frac{(i \alpha)^{j-2}}{j!}
   \sum_{m \le y} \frac{1}{m} \sum_{n \mid m} \frac{\mu(n)}{n}  \\
   & \cdot \sum_{h,k \le \frac{y}{m}} \frac{d_{r}(mh)P([mh]_{y})d_{r}(mk)P([mk]_{y})}{hk}
   \left( \frac{T}{2 \pi h n} \right)^{i \alpha} \frac{nk \rk_{nk}^{(j)}(1)}{\phi(nk)}  \ .
   \label{eq:I22a}
\end{split}
\end{equation}

\subsection{Evaluation of $I_{1}$}
By~(\ref{eq:I1b}) it follows that
\begin{equation}
   I_{1} = \frac{T}{2 \pi}
   \left( -La_{0,0,1} + a_{1,0,1} + a_{0,1,1} - a_{0,0,2} \right)
   + O(TL^{(r+1)^{2}})
   \label{eq:I1c}
\end{equation}
where for $u,v,w \in \mathbb Z_{\ge 0}$ we define $a_{u,v,w}$ to be the sum
$$
   \sum_{mh, mk \le y}
               \frac{\phi(m) d_{r}(mh)P([mh]_{y})d_{r}(mk)P([mk]_{y}) (\log
               h)^{u} (\log k)^{v}}{m^{2}hk} \sum_{p \mid k} p^{i
               \alpha} (\log p)^{w} \ .
$$
By ~(\ref{eq:I1c}) it suffices to evaluate $a_{u,v,w}$. Inverting summation we have
\begin{equation}
\begin{split}
   a_{u,v,w}  =
   \sum_{m \le y} \frac{\phi(m)}{m^{2}}
            \sum_{p \le \frac{y}{m}} \frac{p^{i \alpha} (\log p)^{w}}{p}
    &  \left( \sum_{h \le \frac{y}{m}} \frac{d_{r}(mh)P([mh]_{y})(\log
   h)^{u}}{h} \right)  \\
     &  \cdot \left( \sum_{k \le \frac{y}{pm}} \frac{d_{r}(mpk)P([mpk]_{y})(\log
   pk)^{v}}{k} \right) 
         \nonumber
\end{split}
\end{equation}
By Lemma 10, equations~(\ref{eq:div1}) and~(\ref{eq:div2}) we have
\begin{equation}
\begin{split}
  a_{u,v,w} & = \frac{1}{(r-1)!^2}
   \sum_{mp \le y} \frac{\phi(m)\sigma_{r}(m)^{2}p^{i \alpha} (\log p)^{w}}{m^{2}p}
   \log \left( \frac{y}{m} \right)^{r+u}
   \log \left( \frac{y}{pm} \right)^{r} \\
   & \cdot \int_{0}^{1} \int_{0}^{1} F_{1}(\theta_{1},m) F_{2}(\theta_{2},pm)  \, d \theta_{1}d \theta_{2}
   + \epsilon_{1} + \epsilon_{2} + \epsilon_{3}
   \nonumber
\end{split}
\end{equation}
\begin{equation}
\begin{split}
  & \mathrm{where} \ \epsilon_{1}  \ll \sum_{m\le y}
  \frac{\sigma_{r}(m) L^{u+r}}{m} \sum_{p \le y}
   \frac{(\log
   p)^{w} \epsilon(m)}{p} 
   \ , \\
   & \epsilon_{2}  \ll
   \sum_{m \le y} \frac{\epsilon(m)}{m}
   \sum_{p \le y} \frac{(\log p)^{w}\sigma_{r}(pm)L^{v+r}}{p}  \ ,  \
   \epsilon _{3} \ll
   \sum_{m \le y} \frac{\epsilon(m)}{m}
   \sum_{p \le y} \frac{(\log p)^{w} \epsilon(m)}{p}
   \ .  
   \nonumber
\end{split}
\end{equation}
By~(\ref{eq:sigrbd}) it follows that
\begin{equation}
   \epsilon_{1} \ll L^{u+w+r} \sum_{m \le y} \frac{d_{r}(m)^{2}j_{1}(m)j_{\tau_{0}}(m)}{m}
    \ll L^{u+w+r} \sum_{m \le y} \frac{d_{r}(m)^{2}}{m} \ll L^{u+w+r^2+r} \ .
    \nonumber
\end{equation}
A similar calculation gives $\epsilon_{2} \ll L^{v+w+r^2+r}$ and $\epsilon_{3} \ll L^{w+r^2}$.
Recalling~(\ref{eq:Fi}) and rearranging a little, yields
\begin{equation}
\begin{split}
   & a_{u,v,w} = \frac{(\log y)^{2r+u+v}}{((r-1)!)^2}
    \int_{0}^{1} \int_{0}^{1} \theta_{1}^{r+u-1} \theta_{2}^{r-1}
   \sum_{p \le y} \frac{p^{i \alpha} (\log p)^{w}}{p}   \\
  & \cdot \sum_{m \le \frac{y}{p}} \frac{\phi(m) \sigma_{r}(m) \sigma_{r}(pm)}{m^{2}}
   g_{u,v} \left([m]_{y}, [p]_{y} \right)  d\theta_{1} d\theta_{2}
   + O(L^{\max(u,v)+r^2+r})
   \ .
   \label{eq:auvw2}
\end{split}
\end{equation}
\begin{multline}
  \mathrm{where} \  g_{u,v}(\delta,\beta)   =  \\
   (1-\delta)^{r+u}(1-\beta-\delta)^{r} \left(
  \beta + \theta_{2}(1-\beta-\delta)
  \right)^{v}
  P(\delta+\theta_{1}(1-\delta)) P(\delta + \beta + \theta_{2}(1-\beta-\delta)) \ .
\end{multline}
By Lemma 11$(ii)$,~(\ref{eq:auvw2}) becomes
\begin{equation}
\begin{split}
   & a_{u,v,w} =  rC_{r} (\log y)^{r^2+2r+u+v} \int_{0}^{1} \int_{0}^{1} \theta_{1}^{r+u-1} \theta_{2}^{r-1}
   \sum_{p \le y} \frac{p^{i \alpha}(\log p)^{w}}{p}  \\
  & \cdot \int_{0}^{1-[p]_{y}}  \delta^{r^2-1} g_{u,v} \left(\delta, [p]_{y} \right) \, d\delta \ d \theta_{1}
  d \theta_{2} + \epsilon_{4} + O(L^{\max(u,v)+r^2+r})
   \label{eq:auvwd}
   \nonumber
\end{split}
\end{equation}
where $C_{r}$ is defined by~(\ref{eq:ar}) and
\[
  \epsilon_{4} \ll L^{2r+u+v} \sum_{p \le y} \frac{(\log p)^{w}}{p}
  ( L^{r^{2}}p^{-1} + L^{r^2-1})
  \ll L^{r^2+2r+u+v+w-1}
\]
since $w \ge 1$.  Inverting summation
\begin{equation}
\begin{split}
  & a_{u,v,w} = rC_{r} (\log y)^{r^2+2r+u+v}
  \int_{0}^{1} \int_{0}^{1} \int_{0}^{1} \theta_{1}^{r+u-1} \theta_{2}^{r-1}
  \delta^{r^2-1}  \\
  & \cdot \sum_{p \le y^{1-\delta}} \frac{p^{i \alpha}(\log p)^{w}}{p}
  g_{u,v} \left(\delta, [p]_{y} \right) d\delta d\theta_{1} d \theta_{2}
  + O(L^{r^2+2r+u+v+w-1}) \ .
  \nonumber
\end{split}
\end{equation}
An application of Lemma 12 yields
\begin{equation}
\begin{split}
  & a_{u,v,w} = rC_{r} (\log y)^{r^2+2r+u+v+w}
  \sum_{j=0}^{\infty} \frac{(i \alpha \log y)^{j}}{j!}   \\
  & \cdot \int_{0}^{1} \int_{0}^{1} \int_{0}^{1} \int_{0}^{1-\delta}
  \theta_{1}^{r+u-1} \theta_{2}^{r-1}
  \delta^{r^2-1} \beta^{j+w-1} g_{u,v}(\delta,\beta) d\beta d\delta d \theta_{1} d \theta_{2}
  ( 1 + O(L^{-1})) \ .
  \nonumber
\end{split}
\end{equation}
We write
\[
   \int_{0}^{1} \int_{0}^{1}
   \theta_{1}^{r+u-1} \theta_{2}^{r-1} g_{u,v}(\delta,\beta)  d \theta_{1}  d \theta_{2}
   = (1-\delta)^{r+u}(1-\beta-\delta)^{r} Q_{r+u-1}(\delta) R_{v}(\delta,\beta)
\]
where
\[
  Q_{r+u-1}(\delta) = \int_{0}^{1} \theta_{1}^{r+u-1} P(\delta + \theta_{1}(1-\delta)) \, d \theta_{1} \ ,
\]
\[
  R_{v}(\delta,\beta) = \int_{0}^{1} \theta_{2}^{r-1}
  (\beta+\theta_{2}(1-\delta-\beta))^{v}P(\delta+\beta+\theta_{2}(1-\delta-\beta)) \, d \theta_{2}
\]
and hence
\begin{equation}
\begin{split}
    & a_{u,v,w}  \sim rC_{r}(\log y)^{r^2+2r+u+v+w} \sum_{j=0}^{\infty} \frac{(i \alpha \log y)^{j}}{j!} \\
  &  \cdot \int_{0}^{1}
  \delta^{r^2-1} (1-\delta)^{r+u}
  Q_{u+r-1}(\delta) \int_{0}^{1-\delta} \beta^{j+w-1} (1-\beta -\delta)^{r} R_{v}(\delta,\beta) \, d\beta d\delta
  \ .
  \nonumber
\end{split}
\end{equation}
Now note that
\[
   R_{0}(\delta,\beta) = Q_{r-1}(\delta+\beta)
\ , \ 
  R_{1}(\delta,\beta) = \beta Q_{r-1}(\delta+\beta) + (1-\delta-\beta) Q_{r}(\delta+\beta) \ .
\]
We see that
\begin{equation}
\begin{split}
    & a_{u,0,w}  \sim rC_{r}(\log y)^{r^2+2r+u+w} \sum_{j=0}^{\infty} \frac{(i \alpha \log y)^{j}}{j!} \\
  & \cdot \int_{0}^{1} \int_{0}^{1-\delta}
  \delta^{r^2-1} (1-\delta)^{r+u}
   (1-\beta -\delta)^{r}  \beta^{j+w-1}
   Q_{u+r-1}(\delta) Q_{r-1}(\delta+\beta)
   \, d\beta d\delta
   \nonumber
\end{split}
\end{equation}
and
\begin{equation}
\begin{split}
    & a_{u,1,w}  \sim rC_{r}(\log y)^{r^2+2r+1+u+w} \sum_{j=0}^{\infty} \frac{(i \alpha \log y)^{j}}{j!} \\
  & \cdot \left( \int_{0}^{1} \int_{0}^{1-\delta}
   \delta^{r^2-1} (1-\delta)^{r+u}
  (1-\beta -\delta)^{r} \beta^{j+w}
   Q_{u+r-1}(\delta) Q_{r-1}(\delta+\beta)
   \, d\beta d\delta +  \right. \\
   &  \left.
   \int_{0}^{1} \int_{0}^{1-\delta}
   \delta^{r^2-1} (1-\delta)^{r+u}
  (1-\beta -\delta)^{r+1} \beta^{j+w-1}
   Q_{u+r-1}(\delta) Q_{r}(\delta+\beta)
   \, d\beta d\delta \right) \ .
   \nonumber
\end{split}
\end{equation}
For $\vec{n}=(n_{1},n_{2},n_{3},n_{4},n_{5}) \in \left( \mathbb{Z}_{\ge 0} \right)^{5}$ we 
recall the definition~(\ref{eq:i1234})
\[
  i_{P}(\vec{n}) =
  \int_{0}^{1} \int_{0}^{1-x_{1}} x_{1}^{r^2-1} (1-x_{1})^{n_{1}}(1-x_{1}-x_{2})^{n_{2}}
   x_{2}^{n_{3}} Q_{n_{4}}(x_{1}) Q_{n_{5}}(x_{1}+x_{2}) \, dx_{2} dx_{1}
\]
and hence
\[
  a_{0,0,1} = rC_{r}L^{(r+1)^2}
            \sum_{j=0}^{\infty} \frac{z^{j} \eta^{j+(r+1)^2}}{j!}
            i_{P}(r,r,j,r-1,r-1) \ ,
\]
\[
   a_{1,0,1} = rC_{r} L^{(r+1)^2+1}
             \sum_{j=0}^{\infty} \frac{z^{j} \eta^{j+(r+1)^2+1}}{j!}
             i_{P}(r+1,r,j,r,r-1) \ ,
\]
\[
   a_{0,0,2} = r C_{r} L^{(r+1)^2+1}
             \sum_{j=0}^{\infty} \frac{z^{j} \eta^{j+(r+1)^2+1}}{j!}
             i_{P}(r,r,j+1,r-1,r-1) \ ,
\]
\[
   a_{1,0,1} = r C_{r} L^{(r+1)^2+1}
             \sum_{j=0}^{\infty} \frac{z^{j} \eta^{j+(r+1)^2+1}}{j!}
             i_{P}(r,r,j+1,r-1,r-1) + i_{P}(r,r+1,j,r-1,r) \ .
\]
Combining these identities with~(\ref{eq:I1c}) we arrive at 
\begin{align}
 &  I_{1} \sim rC_{r} \frac{T}{2 \pi} L^{(r+1)^2+1} \sum_{j=0}^{\infty} \frac{z^{j} \eta^{j+(r+1)^2}}{j!}  \label{eq:I1f}\\
 &  \cdot (-i_{P}(r,r,j,r-1,r-1) + \eta (i_{P}(r+1,r,j,r,r-1) + i_{P}(r,r+1,j,r-1,r)) )
 \nonumber
\end{align}
and this is valid up to an error which is smaller by a factor $O(L^{-1})$.

\subsection{Evaluation of $I_{21}$}
We recall that
\begin{equation}
\begin{split}
   I_{21} & \sim \frac{T}{2 \pi}
   \sum_{m \le y} \frac{1}{m} \sum_{n \mid m} \frac{\mu(n)}{n}
   \sum_{h,k \le \frac{y}{m}} \frac{d_{r}(mh)P([mh]_{y})d_{r}(mk)P([mk]_{y})}{hk}  \\
   &
   \cdot \left(
   \frac{1 + i \alpha \log \frac{T}{2 \pi hk n^{2}} -
   \left( \frac{T}{2 \pi h n} \right)^{i \alpha}
   (1-(i \alpha)\log nk)}{(i \alpha)^{2}}
   \right)   \ .
   \nonumber 
\end{split}
\end{equation}
 A little algebra shows that the expression within the brackets simplifies to
\[
   \log \left( \frac{T}{2 \pi h n} \right) \log (nk)
   - (1- (i \alpha) \log nk) \log \left( \frac{T}{2 \pi hn} \right)^{2}
   \sum_{j=0}^{\infty} \frac{(i \alpha \log \left( \frac{T}{2 \pi
   hn} \right))^{j}}{(j+2)!} \ .
\]
We may replace $\log \frac{T}{2 \pi h n}$ by $\log \frac{T}{2 \pi h}$ and $\log (nk)$ by $\log k$ up to an error
of $L(\log n)$. This error term contributes $O(TL^{(r+1)^2})$ as long as we use $|\alpha| \le cL^{-1}$.  It  thus follows that
\begin{equation}
\begin{split}
   I_{21} & \sim  \frac{T}{2 \pi}
   \sum_{m \le y} \frac{\phi(m)}{m^{2}}
   \sum_{h,k \le \frac{y}{m}} \frac{d_{r}(mh)P([mh]_{y})d_{r}(mk)P([mk]_{y})}{hk}  \\
   &
   \cdot \left(
   \log \left( \frac{T}{2 \pi h} \right) \log k
   - (1-(i \alpha) \log k) \log \left( \frac{T}{2 \pi h}
   \right)^{2}
   \sum_{j=0}^{\infty} \frac{(i \alpha \log \left( \frac{T}{2 \pi
   h} \right))^{j}}{(j+2)!} \right) \ .
   \nonumber
\end{split}
\end{equation}
and hence
\begin{equation}
   I_{2} \sim \frac{T}{2 \pi}
   \left(
   b_{1,1}
   -  \sum_{j=0}^{\infty} \frac{(i \alpha)^{j}}{(j+2)!} b_{j+2,0}
   +
   \sum_{j=0}^{\infty} \frac{(i \alpha)^{j+1}}{(j+2)!} b_{j+2,1}
   \right)
   \label{eq:I2buv}
\end{equation}
where
\begin{equation}
    b_{u,v} = \sum_{m \le y} \frac{\phi(m)}{m^{2}}
             \sum_{h,k \le \frac{y}{m}} \frac{d_{r}(mh)P([mh]_{y})(\log
             \frac{T}{2 \pi h})^{u} d_{r}(mk)P([mk]_{y}) 
             (\log k)^{v}}{hk}
    \label{eq:buv}
\end{equation}
for $u,v \ge 0$.  It suffices to evaluate $b_{u,v}$. 
Inserting equations~(\ref{eq:div3}) and~(\ref{eq:div1}) of Lemma 10 in~(\ref{eq:buv}) gives
\begin{multline}
  b_{u,v} \sim \\
  \frac{(\log y)^{u+r}}{((r-1)!)^2}
  \sum_{m \le y} \frac{\phi(m) \sigma_{r}(m)^2}{m^2} \log \left( \frac{y}{m} \right)^{r+v}
  \int_{0}^{1-[m]_{y}} F_{3}(\theta_{1},m) d \theta_{1}
  \int_{0}^{1} F_{1}(\theta_{2},m) d \theta_{2} \
  \nonumber
\end{multline}
where $F_1,F_3$ are given by~(\ref{eq:Fi}).
This is valid up to an error of $O(L^{r^2+r+\max(u,v)})$
and the calculation is analogous to the calculation we did in 
the last section for $a_{u,v,w}$.
Exchanging summation order
and recalling~(\ref{eq:Fi}) gives
\begin{multline}
  b_{u,v} = \\
  \frac{(\log y)^{2r+u+v}}{((r-1)!)^2}
  \int_{0}^{1}\int_{0}^{1} \theta_{1}^{r-1} (\eta^{-1}-\theta_{1})^{u} \theta_{2}^{r+v-1}
  \sum_{m \le y^{1-\theta_{1}}} \frac{\phi(m) \sigma_{r}(m)^{2}}{m^2}
  g \left( [m]_{y} \right)
  d \theta_{1} d \theta_{2}
  \nonumber
\end{multline}
where $g(\delta) = (1-\delta)^{r+v}P(\delta+\theta_{1})P(\delta+(1-\delta)\theta_{2})$. By Lemma 11$(ii)$ we have
\begin{equation}
\begin{split}
  & b_{u,v} = C_{r} (\log y)^{r^2+2r+u+v}
   \int_{0}^{1} \int_{0}^{1}   \theta_{1}^{r-1} (\eta^{-1}-\theta_{1})^{u} \theta_{2}^{r+v-1}
   \int_{0}^{1-\theta_{1}} \delta^{r^2-1} g(\delta)
   d \delta d \theta_{1} d \theta_{2}   
   \nonumber
\end{split}
\end{equation}
plus an error $O(L^{r^2+r+\max(u,v)})$.  Since 
$Q_{r+v-1}(\delta)= \int_{0}^{1} \theta_2^{r+v-1}P(\delta + (1-\delta) \theta_2) \, 
d \theta_2$ it follows that 
\[
  b_{u,v} \sim  C_{r}(\log y)^{r^2+2r+u+v} k_{P}(u,r+v,r+v-1)
\]
where we recall~(\ref{eq:k1234})
\[
  k_{P}(n_{1},n_{2},n_{3})
  = \int_{0}^{1} \int_{0}^{1-\theta_{1}}
    \theta_{1}^{r-1} (\eta^{-1}-\theta_{1})^{n_{1}} \delta^{r^2-1}(1-\delta)^{n_{2}} P(\theta_{1}+\delta)
    Q_{n_{3}}(\delta) d\delta d\theta_{1} \ .
\]
We conclude
\begin{multline}
   I_{21} \sim 
   C_{r} \frac{T}{2 \pi} (\log y)^{r^2+2r+2}
   \sum_{j=0}^{\infty} (z \eta)^{j}
   \left(
   \frac{k_{P}(j+1,r+1,r)}{(j+1)!} - \frac{k_{P}(j+2,r,r-1)}{(j+2)!}
   \right) \ . 
   \label{eq:I21f}
\end{multline}
It can be checked that the error term $O(L^{r^2+r+\max(u,v)})$ contributes
an amount $O(L^{-1})$ smaller than the main term.  

\subsection{Evaluation of $I_{22}$}

By~(\ref{eq:I22a})
\begin{equation}
  I_{22} = - \frac{T}{2 \pi} \sum_{j=2}^{r} \frac{(i \alpha)^{j-2}}{j!} \sum_{u=0}^{\infty}
   \frac{(i \alpha)^{u}}{u!} \, c_{u,j} \label{eq:I22}
\end{equation}
where
\begin{equation}
\begin{split}
   c_{u,j} =  \sum_{m \le y} \frac{1}{m} \sum_{n \mid m}
   \frac{\mu(n)}{n^{i \alpha}}
   & \left(
   \sum_{h \le \frac{y}{m}} \frac{d_{r}(mh)P([mh]_{y})}{h} \log \left( \frac{T}{2 \pi h} \right)^{u}
   \right)  \\
   & \cdot \left(
   \sum_{k\le \frac{y}{m}} d_{r}(mk)P([mk]_{y}) \frac{\rk_{nk}^{(j)}(1)}{\phi(nk)}
   \right) 
    \label{eq:cuj}
\end{split}
\end{equation}
Applying partial summation to~(\ref{eq:drrn}) yields
\begin{align}
   & \sum_{k\le \frac{y}{m}} d_{r}(mk)P([mk]_{y}) \frac{\rk_{nk}^{(j)}(1)}{\phi(nk)}
     \label{eq:drrn2} \\
   & = \frac{\sigma_{r}(m)(-1)^{j} j! \binom{r}{j} \log \left( \frac{y}{m} \right)^{r+j}}{n(r+j-1)!}
   \int_{0}^{1} \theta^{r+j-1} P([m]_{y}+(1-[m]_{y}) \theta) d \theta  +O(E(y)) 
   \label{eq:drrn2}
\end{align}
where $E(y)$ denotes the error term in~(\ref{eq:drrn}).  
We apply Lemma 10$(iii)$ to the first factor in~(\ref{eq:cuj}) and we apply~(\ref{eq:drrn2}) 
to the second factor of~(\ref{eq:cuj}) to obtain
\begin{equation}
\begin{split}
  c_{u,j} & = \frac{(-1)^{j}j! \binom{r}{j}(\log y)^{u+r}}{(r-1)!(r+j-1)!}
  \sum_{m \le y} \sigma_{r}(m)^{2} \log \left( \frac{y}{m} \right)^{r+j}
  \left( \sum_{n \mid m} \frac{\mu(n)}{n^{1+i \alpha}}  \right)  \\
  &  \cdot \int_{0}^{1-[m]_{y}} F_{3}(\theta_{1},m) d \theta_{1}
  \int_{0}^{1} \theta_{2}^{r+j-1} P([m]_{y}+(1-[m]_{y}) \theta_{2}) \, d \theta_{2}  
   \nonumber
\end{split}
\end{equation}
where $F_{3}(\theta_{1},m) = \theta_{1}^{r-1}(\eta^{-1}-\theta_{1})^{u} P([m]_{y}+\theta_{1})$.
Further simplification gives
\begin{equation}
\begin{split}
  c_{u,j} & = \frac{(-1)^{j} j! \binom{r}{j} (\log y)^{2r+u+j}}{(r-1)!(r+j-1)!}
  \int_{0}^{1} \int_{0}^{1} \theta_{1}^{r-1} \theta_{2}^{r+j-1}(\eta^{-1}-\theta_{1})^{u}
  \sum_{m \le y^{1-\theta_{1}}} \frac{\sigma_{r}(m)^2}{m}  \\
 & \cdot \sum_{n \mid m} \frac{\mu(n)}{n^{1+i \alpha}}  
   \left(1 - [m]_{y} \right)^{r+j}
  P\left( \theta_{1} + [m]_{y} \right)
  P \left( [m]_{y} +( 1- [m]_{y}) \theta_{2}  \right)
  d \theta_{1} d \theta_{2} 
  \nonumber
\end{split}
\end{equation}
Now note 
that 
$
      \sum_{n \mid m} n^{-1-i \alpha}= \frac{\phi(m)}{m} + O(|\alpha| 
      \sum_{n \mid m} n^{-1})$.  
Thus we have       
\begin{align}
  &  c_{u,j} = \frac{(-1)^{j} j! \binom{r}{j} (\log y)^{2r+u+j}}{(r-1)!(r+j-1)!}
  \int_{0}^{1} \int_{0}^{1} \theta_{1}^{r-1} \theta_{2}^{r+j-1}(\eta^{-1}-\theta_{1})^{u}  \label{eq:cuj2} \\
  & \cdot \sum_{m \le y^{1-\theta_{1}}} 
  \frac{ \phi(m)\sigma_{r}(m)^2}{m^2}  
   \left(1 - [m]_{y} \right)^{r+j}
  P\left( \theta_{1} + [m]_{y} \right)
  P \left( [m]_{y} +( 1- [m]_{y}) \theta_{2}  \right)
  d \theta_{1} d \theta_{2}
   \nonumber
\end{align}      
plus an  error term of the shape
\begin{equation}
\begin{split}
   & \ll_{r,j,u} |\alpha| (\log y)^{2r+u+j} \sum_{m \le y}
   \frac{\sigma_{r}(m)^2}{m} \sum_{n \mid m} \frac{1}{n}  \\
      & \ll_{r,j,u} |\alpha| L^{2r+u+j}
    \sum_{n \le y} \frac{\sigma_{r}(n)^2}{n^2}
    \sum_{k \le y/m} \frac{\sigma_{r}(k)^2}{k}
    \ll_{r,j,u} L^{r^2+2r+u+j-1} \ . 
\end{split} 
\end{equation} 
Note that we can write down the constant in the $O$ term explicitly
in terms of $r,j,$ and $u$.  Applying Lemma 11 to the inner sum we derive
\begin{equation}
\begin{split}
  c_{u,j} = & \frac{a_{r+1}(-1)^{j} j! \binom{r}{j} (\log y)^{r^2+2r+u+j}}{(r-1)!(r+j-1)!(r^2-1)!}  \\
  &  \cdot \int_{0}^{1} \int_{0}^{1} \int_{0}^{1-\theta_{1}}  \theta_{1}^{r-1} \theta_{2}^{r+j-1}(\eta^{-1}-\theta_{1})^{u}
  \delta^{r^2-1} R(\delta) d \delta d\theta_{1} d \theta_{2} 
  \nonumber
\end{split}
\end{equation}
where $R(\delta) = (1-\delta)^{r+j} P(\theta_{1}+\delta) P(\delta + (1-\delta)\theta_{2})$ and this is valid up to an error of $O_{r,j,u}(L^{r^2+2r+u+j-1})$.  If we recall the definition
$Q_{u}(\delta) = \int_{0}^{1} \theta_{2}^{u} P(\delta + (1-\delta) \theta_{2}) \, d\theta_{2}$
and then execute the integration in the $\theta_{2}$-variable
this becomes
\begin{equation}
\begin{split}
  c_{u,j} & \sim \frac{a_{r+1}(-1)^{j} j! \binom{r}{j} (\log y)^{r^2+2r+u+j}}{(r-1)!(r+j-1)!(r^2-1)!}  \\
  & \cdot
  \int_{0}^{1} \int_{0}^{1-\theta_{1}}  \theta_{1}^{r-1} (\eta^{-1}-\theta_{1})^{u}
  \delta^{r^2-1} (1-\delta)^{r+j} P(\theta_{1}+\delta) Q_{r+j-1}(\delta) d \delta d\theta_{1} \ .
  \nonumber
\end{split}
\end{equation}
Recalling definitions~(\ref{eq:k1234}) and~(\ref{eq:ar}) we have
\begin{equation}
\begin{split}
  c_{u,j} & = \frac{(r-1)!C_{r}(-1)^{j} j! \binom{r}{j} (\log y)^{r^2+2r+u+j}}{(r+j-1)!}
  k_{P}(u,r+j,r+j-1)   \\
  & 
  + O_{r,j,u}(L^{r^2+2r+u+j-1})
  \ . 
  \label{eq:cuj2}
\end{split}
\end{equation}
Combining~(\ref{eq:I22}) and~(\ref{eq:cuj2}) establishes that $I_{22}$ is $-(r-1)!C_{r}\frac{T}{2 \pi}(\log
y)^{r^2+2r+2}$ multiplied by the series
\begin{equation}
\begin{split}
  & \sum_{j=2}^{r} \frac{(-1)^{j} \binom{r}{j} (i \alpha \log y)^{j-2}}{(r+j-1)!} \sum_{u=0}^{\infty}
   \frac{(i \alpha \log y)^{u}}{u!} k_{P}(u,r+j,r+j-1) \\
  = &  \sum_{j=0}^{r-2} \frac{(-1)^{j} \binom{r}{j+2} (i \alpha \log y)^{j}}{(r+j+1)!} \sum_{u=0}^{\infty}
   \frac{(i \alpha \log y)^{u}}{u!} k_{P}(u,r+j+2,r+j+1) \\
  = &  \sum_{j=0}^{r-2} \frac{(-1)^{j} \binom{r}{j+2}}{(r+j+1)!} \sum_{n=j}^{\infty}
   \frac{(\eta z)^{n}}{(n-j)!} k_{P}(n-j,r+j+2,r+j+1) 
   \nonumber
\end{split}
\end{equation}
where we changed $j-2 \to j$ and then made the variable change $n=u+j$ in the inner sum.  Moreover, we can check that the error term $O_{r,j,u}(L^{r^2+2r+u+j-1})$ when substituted in~(\ref{eq:I22}) is smaller than the main term by 
a factor of $O(L^{-1})$.  We now write
$I_{22}=I_{22}^{'} +I_{22}^{''}$ where $I_{22}^{'}$ is the contribution from the $j=0$ term and $I_{22}^{''}$ is the rest.
\begin{equation}
   I_{22}^{'} =   - \frac{(r-1) C_{r} \frac{T}{2 \pi} L^{(r+1)^2+1}}{2(r+1)}
   \sum_{n=0}^{\infty} \frac{z^{n} \eta^{n+ (r+1)^2+1}}{n!} k_{P}(n,r+2,r+1) 
   \label{eq:I22p}
\end{equation}   
\begin{equation}
\begin{split}   
   I_{22}^{''} = & - (r-1)!C_{r}\frac{T}{2 \pi}L^{(r+1)^2+1}
   \sum_{n=1}^{\infty} z^{n} \eta^{n+(r+1)^2+1}  \\
  & \cdot  \sum_{1 \le  j \le \min(n,r-2)} \frac{(-1)^{j} \binom{r}{j+2}}{(n-j)!(r+j+1)!}
  k_{P}(n-j,r+j+2,r+j+1) \ . 
  \label{eq:I22pp}
\end{split}
\end{equation}

\subsection{Evaluating $I$}

We collect our estimates to conclude the evaluation of $I$. 
Since $I= I_1+I_{21}+I_{22}^{'}+I_{22}^{''}$ plus error terms it follows 
from ~(\ref{eq:I1f}),~(\ref{eq:I21f}),~(\ref{eq:I22p}), and~(\ref{eq:I22pp})
that
\begin{equation}
\begin{split}
  I \sim & C_{r} \frac{T}{2 \pi} L^{(r+1)^2+1} \left( \sum_{j=1}^{\infty} z^{j} \eta^{j+(r+1)^2+1}
  \left( \frac{r \hat{i}(r,\eta,j)}{j!}
  + \hat{k}_{1}(r,\eta,j)+\hat{k}_{2}(r,\eta,j) \right) \right) \\
  & +\mathrm{CT}(I)
  \label{eq:Ifinal}
\end{split}
\end{equation}
where $\mathrm{CT}(I)$ denotes the constant term in the above Taylor series, 
\begin{equation}
  \hat{i}(r,\eta,j) = -i_{P}(r,r,j,r-1,r-1) \eta^{-1} + (i_{P}(r+1,r,j,r,r-1)
                           + i_{P}(r,r+1,j,r-1,r-1))  \ ,
   \nonumber
\end{equation}
\begin{equation}
  \hat{k}_{1}(r,\eta,j) = -\frac{k_{P}(j+2,r,r-1)}{(j+2)!}+ \frac{k_{P}(j+1,r+1,r)}{(j+1)!}
                         -\frac{(r-1)k_{P}(j,r+2,r+1)}{2(r+1)j!} \ ,
       \nonumber
\end{equation}
\begin{equation}
  \hat{k}_{2}(r,\eta,j) = -(r-1)! \sum_{u=1}^{\min(j,r-2)} \frac{(-1)^{u} \binom{r}{u+2}}{(j-u)!(r+u+1)!}
  k_{P}(j-u,r+u+2,r+u+1) \ .
   \nonumber
\end{equation}
Next remark that we may conveniently combine $ \hat{k}(r,\eta,j)= \hat{k}_{1}(r,\eta,j)+ \hat{k}_{2}(r,\eta,j)$
to obtain
\begin{equation}
    \hat{k}(r,\eta,j) = -(r-1)! \sum_{u=-2}^{\min(j,r-2)} \frac{(-1)^{u} \binom{r}{u+2}}{(j-u)!(r+u+1)!}
  k_{P}(j-u,r+u+2,r+u+1) \ .
  \label{eq:khat}
\end{equation}
This completes the evaluation of $I$.  

\subsection{The final details}

We now complete the proof of Theorem 2.  In order to abbreviate the following equations
we put 
\begin{equation}
   \theta= C_{r} \frac{T}{2 \pi} L^{(r+1)^2+1}, \ a = \eta^{(r+1)^2-1}, \
    b= \eta^{(r+1)^2}, \ \mathrm{and} \  c= \eta^{(r+1)^2+2} \ . 
\end{equation}
Recall that the discrete moment we are evaluating satisfies
\begin{equation}
   m(H_{r},T;\alpha) = 2 \mathrm{Re}(I)- \overline{J} + O(yT^{1/2+\epsilon}) 
   \ . 
   \label{eq:mlast}
\end{equation}
Moreover, we showed~(\ref{eq:J2}) that $J = \mathrm{CT}(J)(1+O(L^{-1}))$ where 
\begin{equation}
\begin{split}
      \mathrm{CT}(J)  = - \theta&
   \left(  a
   \int_{0}^{1} \alpha^{r^2-1}(1-\alpha)^{2r} Q_{r-1}(\alpha)^2 \, d \alpha
   \right. \\
   & \left. - 2b
     \int_{0}^{1} \alpha^{r^2-1}(1-\alpha)^{2r+1}
      Q_{r-1}(\alpha) Q_{r}(\alpha) \, d \alpha \right) \ . 
      \label{eq:CTJ}
\end{split}
\end{equation}
We shall now combine~(\ref{eq:Ifinal}) and~(\ref{eq:CTJ}) in~(\ref{eq:mlast})
to finish the proof.  In particular we shall now prove that
 $2\mathrm{CT}(I)= \mathrm{CT}(J)$ and hence 
 $\mathrm{CT}(m(H_{r},T,\alpha))=0$.
 This was expected since the constant term in the Taylor series of $\zeta(\rho+\alpha)$ is zero for each $\rho$.  Moreover, the fact that the constant term must be zero provides a consistency check of our calculation.    
We now verify that 
$2 \mathrm{CT}(I)=\mathrm{CT}(J)$.  Recall that $\mathrm{CT}(I)=
\mathrm{CT}(I_{1})+\mathrm{CT}(I_{21})
+\mathrm{CT}(I_{22}^{'})$.  From~(\ref{eq:I1f}) we have 
\begin{equation}
\begin{split}
  & \mathrm{CT}(I_{1}) = \\
  & r \theta
   \left(
   -b \, i_{P}(r,r,0,r-1,r-1)  + c(i_{P}(r+1,r,0,r,r-1)+i_{P}(r,r+1,0,r-1,r))
   \right) \  . 
   \nonumber
\end{split}
\end{equation}
Each of the above integrals has the form 
\begin{equation}
    \int_{0}^{1}
     \int_{0}^{1-x} x^{r^2-1} (1-x)^{u} (1-y-x)^{v} Q_{u-1}(x)Q_{v-1}(x+y) 
     \, dy dx 
     \label{eq:inform}
\end{equation}
for $(u,v)=(r,r),(r+1,r),(r,r+1)$.  Note that we have the identity
\begin{equation}
     (1-x)^{n+1}Q_{n}(x) = \int_{0}^{1-x} \beta^{n} P(x+\beta) \, d \beta 
     \ .
     \label{eq:Qid} 
\end{equation}
One may deduce from~(\ref{eq:Qid}) that 
\[
     \frac{1}{v}(1-x)^{v+1} Q_{v}(x)
     = \int_{0}^{1-x} (1-x-y)^{v} Q_{v-1}(x+y) \, dy  \ . 
\]
and hence 
\[
   (\ref{eq:inform}) = 
   \frac{1}{v}
       \int_{0}^{1}   x^{r^2-1} (1-x)^{u+v+1} Q_{u-1}(x) Q_{v}(x) \, dx \ . 
\]
It follows that
\begin{align}
  \mathrm{CT}(I_{1})  = \theta \cdot  
  \left( -b    \right. & \int_{0}^{1} x^{r^2-1} (1-x)^{2r+1} Q_{r-1}(x)Q_{r}(x) \, dx 
   \\
  + c &  \left( 
  \int_{0}^{1} \right. x^{r^2-1}(1-x)^{2r+2} Q_{r}(x)^2 \, dx  \\
   + 
    \frac{r}{r+1} & \left.  \left. 
    \int_{0}^{1} x^{r^2-1}(1-x)^{2r+2} Q_{r-1}(x) Q_{r+1}(x) \, dx 
    \right) \right) \\
     \nonumber
\end{align}
By~(\ref{eq:I21f}) we have $\mathrm{CT}(I_{21})= \theta \eta^{(r+1)^2+1}
(k_{P}(1,r+1,r)-(1/2)k_{P}(2,r,r-1))$.  Expanding out the factor
$(\eta^{-1}-\theta_1)^2$ in the definition of $k_{P}$ we have 
\begin{equation}
\begin{split}
   &  \theta \eta^{(r+1)^2+1} k_{P}(2,r,r-1) 
   \sim  \theta \eta^{(r+1)^2+1}   \\
   & \cdot \sum_{j=0}^{2} \binom{2}{j} (-1)^{j} \eta^{-(2-j)}
   \int_{0}^{1} \int_{0}^{1-\theta_{1}}
   \theta_{1}^{r-1+j}
   \delta^{r^2-1} (1-\delta)^{r} 
   P(\delta +\theta_{1})Q_{r-1}(\delta) \, d \delta 
   d \theta_{1} \ .
   \nonumber 
\end{split}
\end{equation}
However, by~(\ref{eq:Qid}) this simplifies to 
\begin{equation}
\begin{split}
  \theta \eta^{(r+1)^2+1} k_{P}(2,r,r-1) \sim  \theta
   \left( \right. & 
  a     \int_{0}^{1} 
   \delta^{r^2-1}(1-\delta)^{2r} Q_{r-1}(\delta)^2
   \, d \delta    \\
   - & 2b    \int_{0}^{1} 
   \delta^{r^2-1}(1-\delta)^{2r+1} Q_{r-1}(\delta)
   Q_{r}(\delta) \, d \delta  \\
   & \left. 
  c  \int_{0}^{1} 
   \delta^{r^2-1}(1-\delta)^{2r+2} Q_{r-1}(\delta)
   Q_{r+1}(\delta) \, d \delta  \ . 
  \right)
  \label{eq:b20}
\end{split}
\end{equation} 
Moreover, a similar calculation establishes
\begin{equation}
\begin{split}
   & \theta \eta^{(r+1)^2+1} k_{p}(1,r+1,r) \sim \theta \\
   & \left(
    b \int_{0}^{1} \delta^{r^2-1}(1-\delta)^{2r+1}
    Q_{r-1}(\delta) Q_{r}(\delta) \, d \delta
    -c \int_{0}^{1}
    \delta^{r^2-1} (1-\delta)^{2r+2} Q_{r}(\delta)^2 \, d \delta 
   \right) \ .  
   \label{eq:b11}
\end{split}
\end{equation}
Combining~(\ref{eq:b20}) and~(\ref{eq:b11})  establishes
\begin{equation}
\begin{split}
  \mathrm{CT}(I_{2}) = \theta \left( \right.
 - & \frac{a}{2} \int_{0}^{1} \delta^{r^2-1}(1-\delta)^{2r} Q_{r-1}(\delta)^2  \, d \delta  \\
 +  & 2b \int_{0}^{1} \delta^{r^2-1}(1-\delta)^{2r+1} Q_{r-1}(\delta) Q_{r}(\delta)
 \, d \delta \\
 - & \left. c \int_{0}^{1} \delta^{r^2-1} 
 (1- \delta)^{2r+2} (Q_{r}(\delta)^2 + \frac{1}{2}Q_{r-1}(\delta) Q_{r+1}(\delta)) \, d\delta  
 \right)  
 \nonumber
\end{split}
\end{equation}
In a similar way, it follows from~(\ref{eq:I22p})
\[
  \mathrm{CT}(I_{22}^{'})
  = - \theta \frac{(r-1)}{2(r+1)} c 
  \int_{0}^{1} \delta^{r^2-1}(1-\delta)^{2r+2} Q_{r-1}(\delta)Q_{r+1}(\delta) \, d \delta \ . 
\]
Combining constant terms yields
$\mathrm{CT}(I) = \theta ( c_{1}a +c_{2}b+c_{3}c)$ 
where 
\begin{equation}
\begin{split}
   c_{1} & = - \frac{1}{2} \int_{0}^{1} \delta^{r^2-1}(1-\delta)^{2r} 
   Q_{r-1}(\delta)^2 \, d \delta  \ ,  \\
   c_{2} & = \int_{0}^{1} \delta^{r^2-1}(1-\delta)^{2r+1} Q_{r-1}(\delta)
   Q_{r}(\delta) \, d \delta \ , \\
   c_{3} & = 
   \int_{0}^{1} \delta^{r^2-1} (1-\delta)^{2r+2}  \\
   & \cdot \left( Q_{r}(\delta)^2 (1-1) + Q_{r-1}(\delta) Q_{r+1}(\delta)
    \left( \frac{r}{r+1} - \frac{1}{2}-\frac{(r-1)}{2(r+1)}  \right) \right) \, d \delta  \ . 
     \nonumber
\end{split}
\end{equation}
Observe that $c_{3}=0$ and hence we have shown that
\begin{equation}
\begin{split}
    & \mathrm{CT}(I) = \theta \cdot \\
    & \left(  - \frac{a}{2}
   \int_{0}^{1} \delta^{r^2-1}(1-\delta)^{2r} 
   Q_{r-1}(\delta)^2 \, d \delta 
  + b \int_{0}^{1} \delta^{r^2-1}(1-\delta)^{2r+1} Q_{r-1}(\delta)
   Q_{r}(\delta) \, d \delta  \right) \ .  
   \nonumber
\end{split}
\end{equation}
However, glancing back at~(\ref{eq:CTJ}) we see that $2\mathrm{CT}(I)=
\mathrm{CT}(J)$.  By this fact,~(\ref{eq:Ifinal}),~(\ref{eq:khat}), and
~(\ref{eq:mlast}) we finally deduce   
\[
      m(H_r,T;\alpha) \sim  C_{r} \frac{T}{\pi} L^{(r+1)^2+1} 
      \mathrm{Re} \left( \sum_{j=1}^{\infty} (i \alpha L)^{j} \eta^{j+(r+1)^2+1}
  \left( \frac{r \hat{i}(r,\eta,j)}{j!}
  + \hat{k}(r,\eta,j) \right) \right) \ . 
\]


\begin{thebibliography}{99}

\bibitem{CG} J.B. Conrey and A. Ghosh, "Mean values of the
Riemann zeta function, III", Analytic number theory and
Diophantine problems, Stillwater, OK, 1984, Prog. Math. 70,
Birkh\"auser Boston, Boston 1987, 183--203.

\bibitem{CGG0} J.B. Conrey, A. Ghosh, and S.M. Gonek, "A note on gaps between zeros of
the zeta function", \textit{Bull. London Math. Soc.}, \textbf{16} (1984), 421--424.

\bibitem{CGG} J.B. Conrey, A. Ghosh, and S.M. Gonek,
 "Large gaps between zeros of the zeta-function",
 \textit{Mathematika}, \textbf{33} (1986), 216--238.

\bibitem{CGG2} J.B. Conrey, A. Ghosh, and S.M. Gonek, "Simple
zeros of the Riemann zeta function", \textit{Proc. London Math.
Soc. (3)}, \textbf{76} (1998), 497--522.

\bibitem{E} T. Estermann, "On the representation of a number as
the sum of two products", \textit{Proc. London Math. Soc. (2)},
\textbf{31} (1930),  123--133.

\bibitem{H} R.R. Hall, "A Wirtinger type inequality and the
spacing of the zeros of the Riemann zeta-function", \textit{J.
Number Theory}, \textbf{93} (2002),  235--245.

\bibitem{H2} R.R. Hall, "Generalized Wirtinger
inequalities,random matrix theory, and the zeros of the Riemann
zeta-function", \textit{J. Number Theory}, \textbf{93} (2002),
235--245.

\bibitem{Hu} C.P. Hughes, "Random matrix theory and discrete moments of the Riemann zeta function",
\textit{Journal of Physics A:Math. Gen.}, \textbf{36} (2003) 2907--2917.

\bibitem{Mo} H.L. Montgomery, "The pair correlation of zeros of the zeta function",
Proc. Symp. Pure Math. 24, A.M.S, Providence 1973, 181-193.

\bibitem{MO}  H.L Montgomery and A.M. Odlyzko, "Gaps between zeros of 
the zeta function", Topics in Classical Number Theory, \textit{Colloquia
Math. Soc. Janos Bolyai}, \textbf{34} (Budapest 1981).  

\bibitem{Mu} J. Mueller, "On the difference between consecutive zeros of the 
Riemann zeta function", \textit{J. Number Theory}, \textbf{14}, (1982)
327-331.

\bibitem{Ng} N. Ng, "A note on large and small gaps betweens the zeros
of the Riemann zeta function", in preparation. 

\bibitem{Od} A.M. Odlyzko, "On the distribution of spacings between zeros
of zeta functions", Math. of Comp. \textbf{48} (1987), 273-308.

\bibitem{RS} Z. Rudnick and P. Sarnak, "Zeros of principal L-functions and 
random matrix theory", \textit{Duke Math. J}, \textbf{81} (1996), no.2, 
269-322. 

\bibitem{Se} A. Selberg, Note on a paper of L.G. Sathe, J. of the Indian Math. Soc. B, \textbf{18} (1954), 83--87.

\bibitem{Te}
G. Tenenbaum, {\em Introduction to analytic and probabilistic
number theory}, Cambridge University Press, Cambridge, 1995.

\bibitem{Ti} E.C. Titchmarsh, \textit{The Theory of the Riemann Zeta-Function}
(2nd edition, revised by D.R. Heath-Brown),
Oxford Science Publications, (1986)

\end{thebibliography}
\end{document}